\documentclass[11pt]{article}
\usepackage{amsmath, graphicx, amsfonts,amssymb, calrsfs}
\usepackage{color}

\def\bbR{\mathbb{R}}
\def\r{\rho}
\def\s{\sigma}

\def\vrad{\mathrm{vrad}}

\def\e{\epsilon}
\def\cK{\mathcal{K}}

\def\E{\mathcal{E}}
\def\cF{\mathcal{F}}
\def\o{\omega}
\def\ball{B^n_2}
\def\polar{K^\circ}
\def\l{\lambda}
\def\cS{\mathcal{S}}
\def\cL{\mathbf{L}}
\def\O{\Omega}
\def\bK{\mathbf{K}}
\def\bQ{\mathbf{Q}}
\def\OrliczGmix{G_{\vec{\phi}}^{orlicz}}
\def\OrliczAmix{\O_{\vec{\phi}}^{orlicz}}
\def\OrliczGp{G_{p}^{orlicz}}
\def\OrliczAp{\O_{p}^{orlicz}}
\def\OrliczG{G_{\phi}^{orlicz}}
\def\OrliczA{\O_{\phi}^{orlicz}}
\def\be{\begin{equation}}
\def\ee{\end{equation}}
\def\bea{\begin{eqnarray}}
\def\eea{\end{eqnarray}}
\def\bt{\begin{theorem}}
\def\et{\end{theorem}}
\def\bl{\begin{lemma}}
\def\el{\end{lemma}}
\def\br{\begin{remark}}
\def\er{\end{remark}}
\def\bc{\begin{corollary}}
\def\ec{\end{corollary}}
\def\bd{\begin{definition}}
\def\ed{\end{definition}}
\def\bp{\begin{proposition}}
\def\ep{\end{proposition}}
\newtheorem{theorem}{Theorem}[section]
\newtheorem{lemma}{Lemma}[section]
\newtheorem{remark}{Remark}[section]
\newtheorem{proposition}{Proposition}[section]
\newtheorem{corollary}{Corollary}[section]

\newtheorem{definition}{Definition}[section]
\begin{document}
\title{New Orlicz Affine Isoperimetric Inequalities
}

\author{Deping Ye }
\date{}
\maketitle
\begin{abstract} 
The Orlicz-Brunn-Minkowski theory receives considerable attention recently, and many results in the $L_p$-Brunn-Minkowski theory have been extended to their Orlicz counterparts. The aim of this paper is to develop Orlicz $L_{\phi}$ affine and geominimal surface areas for single convex body as well as for multiple convex bodies, which generalize the $L_p$ (mixed) affine and geominimal surface areas -- fundamental concepts in the $L_p$-Brunn-Minkowski theory. Our extensions are different from the general affine surface areas by Ludwig (in  Adv. Math. 224
(2010)). Moreover, our definitions for Orlicz $L_{\phi}$ affine and geominimal surface areas reveal that these affine invariants are essentially the infimum/supremum of $V_{\phi}(K, L^\circ)$, the Orlicz $\phi$-mixed volume of $K$ and the polar body of $L$, where $L$ runs over all star bodies and all convex bodies, respectively, with volume of $L$ equal to the volume of the unit Euclidean ball $\ball$. Properties for the Orlicz $L_{\phi}$ affine and geominimal surface areas, such as, affine invariance and monotonicity, are proved.  Related Orlicz affine isoperimetric inequalities are also established.

 \end{abstract}\section{Introduction}

The beautiful $L_p$-Brunn-Minkowski theory is developed by the combination of volume and the Firey $p$-sum of convex bodies (i.e., convex compact subsets in $\bbR^n$ with nonempty interiors), and has found many applications such as in geometry. Key tools and objects include, for instance, $L_p$ affine isoperimetric inequalities and $L_p$ affine surface areas. $L_p$ affine isoperimetric inequalities provide upper and/or lower bounds for ($L_p$) affine invariants defined on convex bodies in terms of volume, for example, the $L_p$ affine isoperimetric inequalities for $L_p$ centroid and projection bodies established in \cite{LYZ1} (see also \cite{CG,HaberlFranz2009, LutwakZhang1997}). Lutwak, Yang and Zhang in \cite{LYZ2010a, LYZ2010b} extended such inequalities to their Orlicz counterparts, namely the affine isoperimetric inequalities for Orlicz centroid and projection bodies. These Orlicz affine isoperimetric inequalities initiated the study of the Orlicz-Brunn-Minkowski theory, which involves general (convex) functions and naturally generalizes the $L_p$-Brunn-Minkowski theory. 

The literature for the Orlicz-Brunn-Minkowski theory expands quickly. Important contributions include the study of Logarithmic Minkowski and even Orlicz Minkowski problems \cite{bor2013-2, HaberlLYZ};  the log-Brunn-Minkowski-inequality \cite{bor2012}; stronger versions of the Orlicz-Petty projection inequality \cite{bor2013-1};  Orlicz Busemann-Petty
 centroid inequality  \cite{Chen2011, Li2011, Zhu2012} among others. Recently, in their seminal paper \cite{Gardner2014}, Gardner, Hug and Weil built the foundation and provided a general framework for the Orlicz-Brunn-Minkowski theory. In particular, the Orlicz addition of convex bodies was proposed and the Orlicz mixed volume was defined. Moreover, they proved an Orlicz-Brunn-Minkowski inequality, which extends the classical Brunn-Minkowski inequality. As claimed in \cite{Gardner2014}, ``the Orlicz-Brunn-Minkowski theory is the most general possible based on an addition that retains all the basic geometrical properties enjoyed by the $L_p$-Brunn-Minkowski theory." 

In his ground breaking paper \cite{Lu1}, Lutwak introduced the $L_p$ affine surface areas for $p\geq 1$, which extend the classical
affine surface area by Blaschke \cite{Bl1} in 1923. Note that, for $p\geq 1$ and for (smooth) convex body $K$, the $L_p$ affine surface area of $K$, denoted by $as_p(K)$, has the following integral expression (see Section \ref{section 2} for undefined notation) \be\label{Lp affine surface area} as_p(K)=\int_{S^{n-1}}\left[h_K(u)^{1-p}f_K(u)\right]^{\frac{n}{n+p}}\,d\s(u).\ee This beautiful and convenient integral expression plays fundamental roles in extending the $L_p$ affine surface area from $p\geq 1$ to all $-n\neq p\in \bbR$ (see  e.g.,
\cite{MW2, SW4, SW5}). The $L_p$ affine surface area for $-n\neq p\in \bbR$ is affine invariant (i.e., invariant under all invertible linear transforms with unit absolute value of determinant), and is (essentially) the unique valuation with certain properties such as affine invariance and upper semi-continuity for $p>0$ (see \cite{LudR, LR1} for more precise statements). Important applications of the $L_p$ affine surface areas can be found in, e.g., \cite{A1, A2, Gr2, Jenkinson2012, LSW, Paouris2010, Werner2012a, Werner2012b}. One of the most important results regarding the $L_p$ affine surface area is its related $L_p$ affine isoperimetric inequality (see e.g., \cite{Lu1,
WY2008}): {\it Let convex body $K$ have the origin as its centroid and have the same volume as the unit Euclidean ball $\ball$, then $as_p(K)\leq as_p(\ball)$ for $p>0$ and $as_p(K)\geq as_p(\ball)$ for $-n<p<0$, with equality if and only if $K$ is an origin-symmetric ellipsoid.}

Lutwak also defined the $L_p$ geominimal surface area for $p>1$ in \cite{Lu1}, a concept with similar properties to those 
for the $L_p$ affine surface area. However, one {\it cannot} expect that the $L_p$ geominimal surface area of $K$ has similar integral expression to the formula (\ref{Lp affine surface area}) for its affine counterpart $as_p(K)$. In order to extend the $L_p$ geominimal surface area from $p\geq 1$ to all $-n\neq p\in \bbR$, one turns to {\it observation (a): the $L_p$ affine and geominimal surface areas of $K$ for $p\geq 1$ are (essentially) the infimum of $V_p(K, L^\circ)$, the $p$-mixed volume of $K$ and the polar body of $L$, where $L$ runs over all star bodies and all convex bodies, respectively, with volume of $L$ equal to the volume of the unit Euclidean ball $\ball$} (see \cite{Lu1} and formulas (\ref{Lp-affine-surface-area-Lutwak--1}) and (\ref{Lp-geominimal-surface-area-ye--1}) for more details). In \cite{Ye2014}, the author proved similar formulas for the $L_p$ affine surface areas for all $-n\neq p\in \bbR$, which motivate alternative definitions for the $L_p$ affine surface areas of $K$ and definitions for the $L_p$ geominimal surface areas of $K$ for all $-n\neq p\in \bbR$ (see formulas (\ref{Lp-affine-surface-area-Lutwak--1})-(\ref{Lp-geominimal-surface-area-ye--2})). Contributions for $L_p$ geominimal surface areas including related $L_p$ affine isoperimetric inequalities can be found in, e.g., \cite{Lu1, Petty1974, Petty1985, Schneider2013, Ye2014, Zhou2011}.

Ludwig in \cite{Ludwig2009} defined general affine surface areas involving general convex and concave functions, which are natural extensions of the $L_p$ affine surface area based on formula (\ref{Lp affine surface area}). Affine invariance, the valuation property and affine isoperimetric inequalities for general affine surface areas are established in \cite{Ludwig2009}. The monotone properties of general affine surface areas under the Steiner symmetrization were investigated in \cite{Ye2013} and were used to prove stronger affine isoperimetric inequalities related to general affine surface areas. Again, due to lack of integral expressions for $L_p$ geominimal surface areas, one {\it cannot} follow Ludwig's ideas to define the Orlicz geominimal surface area. 

This paper dedicates to develop the Orlicz $L_{\phi}$ affine and geominimal surface areas. Our motivations are the recent developed Orlicz mixed volume in \cite{Gardner2014} and observation (a). As we can see in Definitions \ref{Orlicz affine surface} and \ref{Orlicz geominimal surface}, our definitions for the Orlicz $L_{\phi}$ affine and geominimal surface areas are consistent with the observation (a): {\it for instance, the Orlicz $L_{\phi}$ affine and geominimal surface areas for $\phi\in \Phi$ are (essentially) the infimum of $V_{\phi}(K, L^\circ)$, the Orlicz $\phi$-mixed volume of $K$ and the polar body of $L$, where $L$ runs over all star bodies and all convex bodies, respectively, with volume of $L$ equal to the volume of $\ball$.}  We prove the affine invariance and monotonicity for the Orlicz $L_{\phi}$ affine and geominimal surface areas. Moreover, we establish the following Orlicz affine isoperimetric inequalities. 

\vskip 2mm \noindent {\bf Theorem \ref{isoperimetric:geominimal}}  {\it Let $K$ be a convex body with centroid at the origin and $B_K$ be the origin-symmetric Euclidean ball with volume equal to the volume of $K$. 

\vskip 2mm \noindent (i) For $\phi\in \Phi$, the following affine isoperimetric inequality holds, with equality if and only if $K$ is an origin-symmetric ellipsoid, 
\begin{eqnarray*}
\OrliczA(K)\leq \OrliczG(K)\leq \OrliczG([B_{K^\circ}]^\circ)= \OrliczA([B_{K^\circ}]^\circ). 
\end{eqnarray*} 
If in addition $\phi$ is increasing, \begin{eqnarray*}
\OrliczA(K)\leq \OrliczG(K)\leq \OrliczG(B_K)= \OrliczA(B_K),
\end{eqnarray*} with equality if $K$ is an origin-symmetric ellipsoid. Moreover, if $\phi$ is strictly increasing, equality holds if and only if $K$ is an origin-symmetric ellipsoid. 

\vskip 2mm \noindent 
(ii) For $\phi\in \Psi$, the following affine isoperimetric inequality holds with equality if $K$ is an origin-symmetric ellipsoid,
\begin{eqnarray*}
\OrliczA(K)\geq \OrliczG(K)\geq \OrliczG(B_K)= \OrliczA(B_K).
\end{eqnarray*}  If in addition $\phi$ is strictly decreasing, equality holds if and only if $K$ is an origin-symmetric ellipsoid. }

In Section \ref{section 5}, we provide definitions for the Orlicz mixed $L_{\phi}$ affine and geominimal surface areas as well as the Orlicz $i$-th mixed $L_{\phi}$ affine and geominimal surface areas. We briefly discuss their properties including affine invariance, Alexander-Fenchel type inequality, and affine isoperimetric inequalities. Basic background and notation in convex geometry are provided in Section \ref{section 2} and readers can read the nice book \cite{Sch} by Schneider for more background in convex geometry.

 
\section{Background and Notation}\label{section 2}
Throughout this paper,  $\cK$  denotes the set of $n$-dimensional {\it convex bodies} in $\bbR^n$: convex compact subsets of
$\bbR^n$ with nonempty interior. Write $\cK_0$ and  $\cK_c$ for the subsets of convex bodies with the origin in their interiors and with centroid at the origin respectively. Denote by $\ball$ the unit Euclidean ball and by $S^{n-1}$ the unit sphere in $\bbR^n$. The volume of $\ball$ is written by $\o_n$. We use $SL(n)$ to denote the group of special linear transformations from $\bbR^n$ to $\bbR^n$. That is, $T\in SL(n)$ means that $T$ is a linear transform with $|det(T)|=1$, where $|det(T)|$ refers to the absolute value of the determinant of $T$. The transpose and the inverse of $T$ will be denoted by $T^{*}$ and $T^{-1}$ respectively. We often write $TK$ for $T(K)$ with $K\in \cK_0$. A convex body $\E\in \cK_0$ is said to be an {\it origin-symmetric ellipsoid} if $\E=rT\ball$ for some $r>0$ and $T\in SL(n)$.

Both support function and radial function can be used to uniquely determine a convex body. The {\it support
function} of $K\in \cK_0$, $h_K: S^{n-1}\rightarrow (0, \infty)$, is defined as
$h_K(u)=\max _{x\in K} \langle x,u \rangle$  where $\langle
\cdot, \cdot\rangle$ denotes the standard inner product  in $\bbR^n$ and induces the usual Euclidian norm $\|\cdot\|$. The {\it radial function}
of $K\in \cK_0$, $\r_K: S^{n-1}\rightarrow (0, \infty)$,  is formulated by \be \label{radial:function-1}\r _K(u)=\max \{\l: \l u\in K\}.\ee  In fact, formula (\ref{radial:function-1}) can be used to define the radial function for star bodies. We use $\cS_0$ to denote the set of star bodies (about the origin) in $\bbR^n$. That is, $L\in \cS_0$ means that the line segment from $0$ to any point $x\in L$ is contained in $L$, and the radial function of $L$, $\r_L(\cdot)$ defined by formula (\ref{radial:function-1}), is continuous and positive.  Clearly, $\cK_0\subset \cS_0$. 

For $K\in \cK_0$, one can define $\polar$, the  {\it polar body} of $K$, by
 $\polar=\{y\in \bbR^n: \langle x,y \rangle \leq 1, \forall x\in
K\}.$  Note that  $\polar\in \cK_0$ and $(K^\circ)^\circ =K$ for $K\in \cK_0$. Moreover, if $K\in \cK_0$, then $\r_K(u){h_{K^\circ}(u)}=1$ holds for all $u\in S^{n-1}$. Write $|K|$ for the volume of $K$ and for the Hausdorff content of its appropriate dimension if $K$ is a general subset of $\bbR^n$. Clearly for all $L\in \cS_0$ and all $K\in \cK_0$, one has   $$|L|=\frac{1}{n}\int_{S^{n-1}}\r_L(u)^n\,d\s(u)\ \ and \ \ |\polar|= \frac{1}{n}\int _{S^{n-1}}
\frac{1} {h _K^n(u)}\,d\s(u),$$ where  $\s$
is the usual spherical measure on $S^{n-1}$.  Denote by $\cK_s\subset \cK_0$ the set of convex bodies with Santal\'{o} point at the origin, i.e., $K\in \cK_s \Leftrightarrow \polar \in \cK_c$. 

 The volume radius of $K$, denoted by $\vrad(K)$, is a way to measure the size of $K\in \cK_0$ in terms of volume. It takes the following form \begin{equation}\label{volume:radius} \vrad(K)=\bigg(\frac{|K|}{|\ball|}\bigg)^{1/n}\Longleftrightarrow |K|^{1/n}=\o_n^{1/n}\vrad(K).\end{equation} 
 Clearly, $\vrad(r\ball)=r$ for all $r>0$, and  for all $T\in SL(n)$,  \begin{equation}\label{affine:invariance:volume:radius} \vrad(TK)=\vrad(K).\end{equation}   It is well known that there is a universal (independent of $K$ and $n$) constant $c>0$, such that, $c\leq \vrad(K)\vrad(K^\circ)\leq 1$ for all $K\in \cK_c$ (or $K\in \cK_s$). The upper bound (i.e., the celebrated Blaschke-Santal\'{o} inequality) is tight, and equality holds if and only if $K$ is an origin-symmetric ellipsoid. The lower bound is the famous Bourgain-Milman inverse Santal\'{o} inequality \cite{BM}. Estimates on the constant $c$ can be found in \cite{GK2, Nazarov2012}. Note that finding the precise minimal value for $\vrad(K)\vrad(K^\circ)$ is still an open problem and is known as the Mahler conjecture:  among all origin-symmetric convex bodies (i.e., $K=-K$), the cube is conjectured to be a minimizer for $\vrad(K)\vrad(K^\circ)$; while the simplex is conjectured to be a minimizer for $\vrad(K)\vrad(K^\circ)$ among all convex bodies $K\in \cK_c$ (or $K\in \cK_s$).

  The {\it Firey
$p$-sum} \cite{Firey1962} of  $K, L\in \cK_0$ with  $\l, \eta \geq 0$ (not both zeros) for $p\geq 1$, denoted by $\l K+_p \eta L$, is determined by the support function 
$$\big(h_{\l K+_p\eta L}(u)\big)^p=\lambda
\big(h_K(u)\big)^p+\eta \big(h_L(u)\big)^p.$$ 
Lutwak \cite{Lut1993} defined 
$V_p(K,L)$, the {\it $p$-mixed volume}  of $K$ and $L$,  by  \begin{equation*}
\frac{1}{p}\ V_p(K,L)=\lim _{\e
\rightarrow 0} \frac{|K+_p\e L|-|K|}{n\e}, \end{equation*} and  proved that there is a measure $S_p(K,\cdot)$,
such that
\begin{equation*} V_p(K,L)=\frac{1}{n}\int _{S^{n-1}}h_L^p(u)\,dS_p(K,u).  \end{equation*} Note that $V_1
(K,L)$ for $p=1$ is the classical  {\it mixed volume} of $K$ and $L$. It has been proved that 
there is a positive Borel measure $S(K, \cdot)$ on $S^{n-1}$ (see \cite{Ale1937-1, Fenchel1938}), such
that, for $K, L\in \cK$, $$V_1 (K,L)=\frac{1}{n}\int
_{S^{n-1}}h_L(u)\,dS(K,u).$$ The relation between $S(K, \cdot)$ and $S_p(K,\cdot)$ for all $p\geq 1$ is
\be\,dS_p(K,u)=h_K^{1-p}(u)\,dS(K, u). \nonumber \ee  One can actually define  the $p$-mixed volume (see \cite{Ye2014}) for all $p\in \bbR$ by 
\be V_p(K, L) =\frac{1}{n} \int
_{S^{n-1}} h_L(u)^ph_K^{1-p}(u)\,dS(K, u).\label{p:mixed:general p}\ee When $K\in \cK_0$ and $L\in \cS_0$, we use the following formula    $$ V_p(K, L^\circ)=\frac{1}{n} \int
_{S^{n-1}} \r_L(u)^{-p} h_K^{1-p}(u)\,dS(K, u), \ \ p\in \bbR.$$  
It is easy to see that for all $\lambda>0$, \begin{equation} \label{homogeneous:pmixed volume} V_p(K, \lambda L)=\lambda^p V_p(K, L).\end{equation}  
We say that  $K\in \cK$ has  {\it curvature function} $f_K(\cdot)$ if the measure $S(K, \cdot)$ is absolutely continuous with respect to the spherical measure $\s$ and satisfies $dS(K,u)=f_K(u)d\s(u).$
Let $\cF_0\subset \cK_0$ be the subset of convex bodies with curvature function and with the origin in its interior. Also, let   $\cF_c=\cF_0\cap \cK_c$ and $\cF_s=\cF_0\cap \cK_s$. The set $\cF_0^+$ denotes the subset of convex bodies in $\cF_0$ with continuous positive curvature functions.

 The $L_p$ affine surface area of $K$, $as_p(K)$, is a fundamental object in affine convex geometry. It can be formulated by (see, e.g., \cite{Lu1, Ye2014})    
 \begin{eqnarray}\label{Lp-affine-surface-area-Lutwak--1}
as_p(K)&=&\inf_{L\in \cS_0} \left\{n V_p(K, L^\circ) ^{\frac{n}{n+p}}\
|L|^{\frac{p}{n+p}}\right\}, \ \ \ p\geq 0; \\\label{Lp-affine-surface-area-Lutwak--2}
as_p(K)&=&\sup_{L\in \cS_0} \left\{n V_p(K, L^\circ) ^{\frac{n}{n+p}}\
|L|^{\frac{p}{n+p}}\right\}, \ \ \ -n\neq p<0.
\end{eqnarray}  The $L_p$ affine surface area has many nice properties such as the valuation property. Moreover, for all invertible linear transforms $T: \bbR^n\rightarrow \bbR^n$, one has $$as_p(TK)=|det(T)|^{\frac{n-p}{n+p}}as_p(K). $$  The $L_p$ geominimal surface area of $K$,  $\tilde{G}_p(K)$, can be defined by \cite{Ye2014} \begin{eqnarray}\label{Lp-geominimal-surface-area-ye--1}
\tilde{G}_p(K)&=&\inf_{Q\in \cK_0} \left\{n V_p(K, Q^\circ) ^{\frac{n}{n+p}}\
|Q|^{\frac{p}{n+p}}\right\}, \ \ \ p\geq 0; \\\label{Lp-geominimal-surface-area-ye--2}
\tilde{G}_p(K)&=&\sup_{Q\in \cK_0} \left\{n V_p(K, Q^\circ) ^{\frac{n}{n+p}}\
|Q|^{\frac{p}{n+p}}\right\}, \ \ \ -n\neq p<0.
\end{eqnarray}  Note that for all invertible linear transforms $T: \bbR^n\rightarrow \bbR^n$, one also has $$ \tilde{G}_p(TK)=|det(T)|^{\frac{n-p}{n+p}}\tilde{G}_p(K). $$  Recall that Lutwak defined the $L_p$ geominimal surface area of $K$ for $p\geq 1$ in \cite{Lu1} by $$\o_n^{p/n}G_p(K)=\inf_{Q\in \cK_0}\{n V_p(K, Q)|Q^\circ|^{p/n} \}.$$ Thus, the relation between $\tilde{G}_p(K)$ and  $G_p(K)$ for $p\geq 1$ is  \be [\tilde{G}_p(K)]^{n+p}=(n\o_n)^p [G_p(K)]^n.  \label{relation tilde G}\ee In fact, one has the following formulas for the $L_p$ geominimal surface areas: \begin{eqnarray}\label{Lp-geominimal-surface-area-ye--1--1}
\big(\tilde{G}_p(K)\big)^{\frac{n+p}{n}} &=&\inf_{Q\in \cK_0} \left\{n^{\frac{n+p}{n}}  V_p(K, Q^\circ) \
|Q|^{\frac{p}{n}}\right\}, \ \ \ p\geq 0\ \& \ p<-n; \\ \label{Lp-geominimal-surface-area-ye--2--2}
\big(\tilde{G}_p(K)\big)^{\frac{n+p}{n}} &=&\sup_{Q\in \cK_0} \left\{n^{\frac{n+p}{n}}  V_p(K, Q^\circ) \
|Q|^{\frac{p}{n}}\right\}, \ \ \ -n<p<0.  \end{eqnarray} Similar formulas for $L_p$ affine surface areas can also be obtained. Based on formulas (\ref{Lp-geominimal-surface-area-ye--1--1}) and (\ref{Lp-geominimal-surface-area-ye--2--2}), we have the following {\it observation (b): the convexity and concavity of function $t^{-p/n}$ (not $t^p$ itself) determines the supremum and the infimum.} More precisely, formula (\ref{Lp-geominimal-surface-area-ye--1--1}) is related to $t^{-p/n}$ being convex for $p\geq 0$ and $p<-n$; while formula (\ref{Lp-geominimal-surface-area-ye--2--2}) is related to $t^{-p/n}$ being concave for  $-n<p<0$.

We will work on function $\phi(t): (0, \infty)\rightarrow (0, \infty)$. Such a function is said to be convex if $$\phi(\lambda t+(1-\lambda)s)\leq \lambda \phi(t)+(1-\lambda)\phi(s),  \ \ \forall t, s\in (0, \infty)\ and\ \lambda\in [0,1].$$ Function $\phi(t)$ is said to be strictly convex if $$\phi(\lambda t+(1-\lambda)s)< \lambda \phi(t)+(1-\lambda)\phi(s),  \ \ \forall t, s\in (0, \infty)\ with\ t\neq s, and\ \lambda\in (0,1).$$ Similarly, one can define concave and strictly concave function $\phi(t)$ by changing the directions of above inequalities. Function $\phi$ is increasing if $\phi(t)\leq \phi(s)$ and is strictly increasing if $\phi(t)<\phi(s)$ for all $t<s$; while function $\phi$ is decreasing if $\phi(t)\geq \phi(s)$ and is strictly decreasing if $\phi(t)>\phi(s)$ for all $t<s$. The inverse function of $\phi$, if exists,  is denoted by $\phi^{-1}(t)$.

\section{Orlicz $L_{\phi}$ affine and geominimal surface areas}\label{section 3}

Let $\phi: (0, \infty)\rightarrow (0, \infty)$ be a positive continuous function. Define the Orlicz $\phi$-mixed volume $V_{\phi}(K, Q)$ of convex bodies $K, Q\in \cK_0$ by  $$V_{\phi}(K, Q)=\frac{1}{n}\int _{S^{n-1}}\phi\left(\frac{h_Q(u)}{h_K(u)}\right) h_K(u)\,dS(K, u).$$  When $K\in \cK_0$ and $L\in \cS_0$, we use the following  formula  $$ V_{\phi}(K, L^\circ)=\frac{1}{n} \int
 _{S^{n-1}}\phi\left(\frac{1}{\r_L(u) h_K(u)}\right) h_K(u)\,dS(K, u).$$  For $\phi(t)=t^p$ with $p\in \bbR$, one gets the $p$-mixed volume $V_p(K, Q)$ given by formula (\ref{p:mixed:general p}). If $\phi$ is a convex function satisfies $\phi(0)=0$ and $\phi(1)=1$, our Orlicz $\phi$-mixed volume $V_{\phi}(K, Q)$ is identical to the one introduced in \cite{Gardner2014}. 
 
For function $\phi: (0, \infty)\rightarrow (0, \infty)$, define $F_{\phi}(t)=\phi(t^{-1/n})$. For simplicity, we often write $F(t)$ for $F_{\phi}(t)$ and clearly $\phi(t)=F(t^{-n})$. From observation (b), one sees that it is the convexity and concavity of function $F(t)$ (not $\phi(t)$ itself) determining the supremum and the infimum. This observation is further strengthened by the proof of Corollary  \ref{Orlicz:ellpsoids}. Clearly, if $\phi(t)$ is increasing then $F(t)$ is decreasing; while if $\phi(t)$ is decreasing then $F(t)$ is increasing. Vice verse, if $F(t)$ is increasing then $\phi(t)$ is decreasing; while if $F(t)$ is decreasing then $\phi(t)$ is increasing.  The relation of convexity and concavity between $\phi(t)$ and $F(t)$ is not clear. However, if $\phi(t)$ is convex and increasing, then $F(t)$ is convex and decreasing. To see this, let $t, s\in (0, \infty)$ and $\lambda \in [0, 1]$, then \begin{eqnarray*} F\big(\lambda t+(1-\lambda)s\big)&=&\phi\big([\lambda t+(1-\lambda)s]^{-1/n}\big)\leq \phi\big(\lambda t^{-1/n}+(1-\lambda)s^{-1/n}\big)\\ &\leq& \lambda\phi\big(t^{-1/n}\big)+(1-\lambda)\phi\big(s^{-1/n}\big)= \lambda F(t)+(1-\lambda)F(s),\end{eqnarray*} where the first inequality follows from the convexity of $t^{-1/n}$ and the increasing property of $\phi(t)$, and the second inequality follows from the convexity of $\phi(t)$. Similarly, if $F(t)$ is convex and increasing, then $\phi(t)$ is convex and decreasing.
 
 Let the set of functions $\Phi$ be \begin{eqnarray*} \Phi=\{\mbox{$\phi: (0, \infty)\rightarrow (0, \infty)$ :  $F(t)$ is either a constant or a strictly convex function}\}.\end{eqnarray*}  The set $\Phi$ contains functions such as $t^p$ with $p\in (-\infty, -n)\cup (0, \infty)$ and all convex increasing functions. Note that $\Phi$ could have neither convex nor concave functions, such as,  $e^{-(t^{-n})}$ whose second order derivative is  $e^{-(t^{-n})}nt^{-n-2}(nt^{-n}-n-1)$ and changes the sign at (the inflection point) $t=\big(\frac{n}{n+1}\big)^{1/n}$. One may define even more general set than the above set $\Phi$ by changing ``strictly convex" to convex;  however, the strict convexity of $F(t)$ required in the definition of $\Phi$ is mainly to exclude those functions $\phi$ proportional to $t^{-n}$ or even part of $\phi(t)$ proportional to $t^{-n}$, a function with problems in defining related affine and geominimal surface areas. (Note that $L_p$ affine and geominimal surface areas were defined for functions $\phi(t)=t^p$ with $p\neq -n$, see formulas (\ref{Lp-affine-surface-area-Lutwak--1}), (\ref{Lp-affine-surface-area-Lutwak--2}), (\ref{Lp-geominimal-surface-area-ye--1}) and (\ref{Lp-geominimal-surface-area-ye--2})).
  
Note that if a concave function $F: (0, \infty)\rightarrow (0, \infty)$ is decreasing on $[t_0, \infty)$ for some $t_0>0$, then $F(t)=F(t_0)$ for all $t>t_0$. To this end, one first gets that $\lim_{t\rightarrow \infty} F(t)$ exists and is finite, say equal to $a\geq 0$. On the other hand, let $\lambda=\frac{t-2t_0}{2t-2t_0}$ for $t>2t_0$, one has $ F(t/2)=F(\lambda t+(1-\lambda)t_0)\geq \lambda F(t)+(1-\lambda)F(t_0).$  Taking limit with $t\rightarrow\infty$, one gets $\lambda\rightarrow 1/2$ and hence  $a\geq a/2+F(t_0)/2\Leftrightarrow a\geq  F(t_0). $  This further implies that $a=F(t_0)$ as $F(t)$ is decreasing, which leads to $F(t)=a$ for all $t\geq t_0$. In view of this, let the set of functions $\Psi$ be
\begin{eqnarray*} \Psi=\{\phi: (0, \infty)\rightarrow (0, \infty): F(t)\ \mbox{is either a constant or an increasing strictly concave function}\}.\end{eqnarray*} Clearly, $\phi$ is decreasing and so is (if exists) $\phi^{-1}(t)$ for $\phi\in\Psi$. Sample (non-constant) functions in $\Psi$ are: $t^p$ with $p\in (-n, 0)$,  $\arctan(t^{-n})$, and $\ln (1+t^{-n})$. Note that $\Psi$ could have neither convex nor concave functions, such as,  $\arctan(t^{-n})$ whose second order derivative is $ {nt^{-n-2}(1+t^{-2n})^{-2} [(n+1)-(n-1)t^{-2n}]}$ and changes the sign at (the inflection point) $t=\big(\frac{n-1}{n+1}\big)^{\frac{1}{2n}}$. One may define even more general set than the above set $\Psi$ by changing ``strictly concave" to concave;  however, the strict concavity of $F(t)$ required in the definition of $\Psi$ is mainly to exclude those functions $\phi$ proportional to $t^{-n}$ or even part of $\phi(t)$ proportional to $t^{-n}$.  We also mention that both $\Phi$ and $\Psi$ may not contain some nice functions such as $\phi(t)=e^{-t}$. 

We now propose the following definition for the Orlicz $L_{\phi}$ affine surface area.  \bd \label{Orlicz affine surface}  Let $K\in \cK_0$ be a convex body with the origin in its interior.  \\
 (i)  For $\phi\in \Phi$, we define the Orlicz $L_\phi$ affine surface area of $K$ by \begin{eqnarray*}
  \OrliczA(K)=\inf_{L\in \cS_0} \left\{nV_{\phi}(K, \vrad(L)L^\circ) \right\}=\inf \left\{nV_{\phi}(K, L^\circ):  L\in \cS_0 \ with\ |L|=\o_n \right\}. 
  \end{eqnarray*}  (ii)  For $\phi\in \Psi$, we define the Orlicz $L_\phi$ affine surface area of $K$ by \begin{eqnarray*}
    \OrliczA(K)=\sup_{L\in \cS_0} \left\{nV_{\phi}(K, \vrad(L)L^\circ) \right\}=\sup \left\{nV_{\phi}(K, L^\circ):  L\in \cS_0 \ with\ |L|=\o_n \right\}.  \end{eqnarray*} 
\ed   \noindent {\bf Remark.}  Clearly, if $\phi(t)=a>0$ is a constant function, then $\OrliczA(K)=a n|K|$. We write $\O_{p}^{orlicz}(K)$ for the case $\phi(t)=t^p$ with $-n\neq p\in \bbR$ and in fact $$(n\o_n)^{p/n}\OrliczAp(K)= \big(as_p(K)\big)^{\frac{n+p}{n}},$$ where $as_p(K)$ is the $L_p$ affine surface area of $K$ defined in (\ref{Lp-affine-surface-area-Lutwak--1}) and (\ref{Lp-affine-surface-area-Lutwak--2}).  As an example, we show this by letting $\phi(t)=t^p$ with $p\geq 0$. Then 
  \begin{eqnarray*}
 (n\o_n)^{p/n}\OrliczAp(K)&=&(n\o_n)^{p/n}\inf_{L\in \cS_0} \left\{nV_{p}(K, \vrad(L) L^\circ) \right\}\\&=&  n^{\frac{n+p}{n}}\inf_{L\in \cS_0} \left\{V_{p}(K,  L^\circ) [\o_n^{1/n}\vrad(L)]^{p}\right\}\\&=& \left(\inf_{L\in \cS_0} \left\{nV_{p}(K,  L^\circ)^{\frac{n}{n+p}} |L|^{\frac{p}{n+p}}\right\}\right)^{\frac{n+p}{n}}=\big(as_p(K)\big)^{\frac{n+p}{n}},
  \end{eqnarray*}  where we have used formula (\ref{homogeneous:pmixed volume}) in the second equality and formula (\ref{volume:radius}) in the third equality.

 We propose the following definition for the Orlicz $L_{\phi}$ geominimal surface area. 
\bd \label{Orlicz geominimal surface}  Let $K\in \cK_0$ be a convex body with the origin in its interior. \\
 (i)  For $\phi\in \Phi$, we define the Orlicz $L_\phi$ geominimal surface area of $K$ by
\begin{eqnarray}
\OrliczG(K)=\inf_{Q\in \cK_0} \left\{nV_{\phi}(K, \vrad(Q^\circ)Q) \right\}=\inf \left\{nV_{\phi}(K, Q):  Q\in \cK_0 \ with\ |Q^\circ|=\o_n \right\}.\label{Orlicz-geominimal-surface-area-convex}
\end{eqnarray}
(ii)  For $\phi\in \Psi$, we define the Orlicz $L_\phi$ geominimal surface area of $K$ by 
\begin{eqnarray}
 \OrliczG(K)=\sup_{Q\in \cK_0} \left\{nV_{\phi}(K, \vrad(Q^\circ)Q) \right\}=\sup \left\{nV_{\phi}(K, Q):  Q\in \cK_0 \ with\ |Q^\circ|=\o_n \right\}.\label{Orlicz-geominimal-surface-area:concave}   \end{eqnarray} \ed  
\noindent {\bf Remark.}  Clearly, if $\phi(t)=a>0$ is a constant function, then $\OrliczG(K)=a n|K|$. We write $G_{p}^{orlicz}(K)$ for the case $\phi(t)=t^p$ with $-n\neq p\in \bbR$ and in fact, similar to the remark after Definition \ref{Orlicz affine surface},  $$ (n\o_n)^{p/n}\OrliczGp(K)= \big(\tilde{G}_p(K)\big)^{\frac{n+p}{n}},
 $$ where $\tilde{G}_p(K)$ is given by  formulas (\ref{Lp-geominimal-surface-area-ye--1}) and (\ref{Lp-geominimal-surface-area-ye--2}). Formula (\ref{relation tilde G}) further implies $\OrliczGp (K)= G_p(K)$ for $p\geq 1$, i.e., our Orlicz $L_{\phi}$ geominimal surface area for $\phi(t)=t^p$ with $p\geq 1$ is identical to the one by Lutwak. 
 \bp Let $K\in \cK_0$. Assume that  $\phi\leq \psi$ with either $\phi, \psi\in \Phi$   or  $\phi, \psi\in \Psi$, then \begin{eqnarray*}
\OrliczA(K)\leq  \O_{\psi}^{orlicz} (K)  \ \ and \ \  
\OrliczG(K)\leq  G_{\psi}^{orlicz} (K).
\end{eqnarray*} \ep \noindent {\bf Proof.} We only prove the case for the Orlicz $L_{\phi}$ geominimal surface area with $\phi, \psi\in \Phi$ such that $\phi\leq \psi$. Other cases can be proved along the same line. 

In fact, for all $K, Q\in \cK_0$, one has  $$\phi\left(\frac{h_Q(u)}{h_K(u)}\right)\leq \psi\left(\frac{h_Q(u)}{h_K(u)}\right), \ \ \ \forall u\in S^{n-1}.$$ Hence, $V_{\phi}(K, Q)\leq V_{\psi}(K, Q)$ for all $K, Q\in \cK_0$. Formula (\ref{Orlicz-geominimal-surface-area-convex})  implies that  $\OrliczG(K)\leq  G_{\psi}^{orlicz} (K)$ for $\phi\leq \psi$ with $\phi, \psi\in \Phi$.

The following proposition states that both Orlicz $L_{\phi}$ affine and geominimal surface areas are affine invariant. 
 \bp\label{homogeneous:degree}  Let $K\in \cK_0$. For all $\phi\in \Phi$ or $\phi\in \Psi$, one has $$ \OrliczA(TK)=\OrliczA(K); \ \ \ \OrliczG(TK)=\OrliczG(K), \ \  \forall T\in SL(n).$$ \ep 
\noindent {\bf Proof.} Let $T\in SL(n)$ and $u=u(v)=\frac{T^{*} v}{\|T^{*} v\|}\in S^{n-1}$ for $v\in S^{n-1}$. Note that $\frac{1}{n} h_K(u)\,dS(K, u)$ is the volume element of $K$ and hence  $ h_{TK}(v)\,dS(TK, v)=h_K(u)\,dS(K, u)$. On the other hand, $\forall y\in TK$, there is a unique $x\in K$ s.t. $y=Tx$ (as $T$ is invertible). Thus,\begin{equation}\label{support:T:function} h_{TK}(v)=\max_{y\in TK} \langle y, v\rangle=\max_{x\in K} \langle Tx, v\rangle =\max_{x\in K} \langle x, T^{*} v\rangle=\|T^*v\| \max_{x\in K} \langle x, u \rangle=\|T^*v\| h_K(u),\end{equation} which implies that,  for $ T\in SL(n)$, 
\begin{eqnarray*}
nV_{\phi}(TK, TQ)\!=\!\! \int _{S^{n-1}} \!\phi \bigg(\!\frac{h_{TQ}(v)}{h_{TK}(v)}\!\bigg) h_{TK}(v)\,dS(TK, v)\!=  \!\! \int _{S^{n-1}}\! \phi\bigg(\!\frac{h_{Q}(u)}{h_{K}(u)}\!\bigg) h_{K}(u)\,dS(K, u) \!=\! n  V_{\phi}(K, Q). \end{eqnarray*} Together with formula (\ref{affine:invariance:volume:radius}), one gets, for all $\phi\in \Phi$, \begin{eqnarray*}
\OrliczG(TK)=\inf_{Q\in \cK_0} \left\{nV_{\phi}\big(TK, \vrad((TQ)^\circ)(TQ)\big) \right\} =\inf_{Q\in \cK_0} \left\{nV_{\phi}(K, \vrad( Q ^\circ) Q) \right\}=\OrliczG (K);
\end{eqnarray*} while for all $\phi\in \Psi$, \begin{eqnarray*}
\OrliczG(TK)=\sup_{Q\in \cK_0} \left\{nV_{\phi}\big(TK, \vrad((TQ)^\circ)(TQ)\big) \right\} =\sup_{Q\in \cK_0} \left\{nV_{\phi}(K, \vrad( Q ^\circ) Q) \right\}=\OrliczG (K).
\end{eqnarray*}
The proof for the case $\OrliczA(TK)=\OrliczA(K)$ follows along the same line if one notices that, by similar calculation to formula (\ref{support:T:function}),  $\r_{(T^*)^{-1} L}(v){\|T^*v\|}= {\r_L(u)}$ and therefore
\begin{eqnarray*}
nV_{\phi}(TK, ((T^*)^{-1}L)^\circ )&=& \int _{S^{n-1}}  \phi \bigg( \frac{1}{\r_{(T^*)^{-1} L}(v)h_{TK}(v)}\!\bigg) h_{TK}(v)\,dS(TK, v)\\ &=&\int _{S^{n-1}}\! \phi\bigg( \frac{1}{\r_L(u)h_{K}(u)} \bigg) h_{K}(u)\,dS(K, u)  =  n  V_{\phi}(K, L^\circ). \end{eqnarray*} 

\noindent {\bf Remark.} When $\phi(t)=t^p$ with $-n\neq p\in \bbR$, a more careful calculation shows that  $$ \OrliczGp (TK)=|det(T)|^{\frac{n-p}{n}} \OrliczGp(K); \ \ \OrliczAp (TK)=|det(T)|^{\frac{n-p}{n}} \OrliczAp(K).$$  Moreover, if $p=1$, both $ \OrliczAp(K)$ and $ \OrliczGp(K)$ are translation invariant, i.e., $\forall z_0\in \bbR^n$, $$\OrliczGp (K-z_0)=\OrliczGp(K); \ \ \OrliczAp (K-z_0)=\OrliczAp(K), \ for \  p=1.$$ However, for general $\phi\in \Phi$ or $\phi\in \Psi$, one {\it cannot} expect the translation invariance even for $K$ being ellipsoids, see inequalities (\ref{non translation invariance 1}) and (\ref{non translation invariance 2}) for special cases.  
\bp \label{Comparison:geominimal:affine} Let $K\in \cK_0$ be a convex body in $\bbR^n$ with the origin in its interior. \\ (i) For $\phi\in \Phi$, one has 
$$\OrliczA(K)\leq \OrliczG(K).$$
(ii) For $\phi\in \Psi$, one has 
$$\OrliczA(K)\geq \OrliczG(K).$$
 \ep  \noindent {\bf Proof.}  (i). Note that $\cK_0\subset\cS_0$. Therefore, for $\phi\in \Phi$, one has 
  \begin{eqnarray*} 
   \OrliczA(K)=\inf_{L\in \cS_0} \left\{nV_{\phi}(K, \vrad(L) L^\circ) \right\} \leq   \inf_{Q\in \cK_0} \left\{nV_{\phi}(K, \vrad(Q) Q^\circ) \right\}= \OrliczG(K).
   \end{eqnarray*}   (ii). For $\phi\in \Psi$, one gets \begin{eqnarray*} 
   \OrliczA(K)=\sup_{L\in \cS_0} \left\{nV_{\phi}(K, \vrad(L) L^\circ) \right\} \geq   \sup_{Q\in \cK_0} \left\{nV_{\phi}(K, \vrad(Q) Q^\circ) \right\}=  \OrliczG(K).
   \end{eqnarray*}
 
 \bc\label{Orlicz:ellpsoids} Let $\E$ be an origin-symmetric ellipsoid. For $\phi\in\Phi$ or $\phi\in\Psi$, one has   $$\OrliczA(\E)=\OrliczG(\E)=n\phi\left(  {\vrad(\E^\circ)} \right) |\E|. $$\ec
 \noindent {\bf Proof.}  Let us first calculate $\OrliczA(r\ball)$ for some $r>0$. Taking $Q=\ball$ (and hence $\vrad(\ball)=1$) in formula (\ref{Orlicz-geominimal-surface-area-convex}), one has, for $\phi\in \Phi$, \begin{eqnarray*}
  \OrliczA(r\ball) \leq  \OrliczG(r\ball) \leq  n  V_{\phi}(r\ball ,  \ball)  =  \int _{S^{n-1}}\phi\left(1/r\right) r^n\,d\s(u)=n \phi (1/r) |r\ball|.
  \end{eqnarray*}On the other hand, Proposition \ref{Comparison:geominimal:affine} implies   \begin{eqnarray*}
    \OrliczG(r\ball)&\geq& \OrliczA(r\ball)= n \inf_{L\in \cS_0}V_{\phi}(r\ball , \vrad(L) L^\circ)  \\   &=&  nr^n \inf_{L\in \cS_0}\left[ \frac{1}{n}\int _{S^{n-1}}F\left(\frac{r^n\r_L^n(u)|\ball|}{|L|}\right)\,d\s(u) \right]\\ &\geq & n|r\ball| \inf_{L\in \cS_0}F\left( r^n\cdot \frac{1}{n|L|}\int _{S^{n-1}} {\r_L^n(u)}  \,d\s(u) \right)\\ &=& nF\left( r^n\right) |r\ball|  =  n \phi\left( 1/r \right) |r\ball|,
    \end{eqnarray*} where the second inequality follows from Jensen's inequality (for convex function $F$ as $\phi\in \Phi$). In conclusion, for $\phi\in\Phi$,  \begin{eqnarray*}
     \OrliczG(r\ball)= \OrliczA(r\ball)=  n  \phi\left(1/r \right) |r\ball|=n\phi (\vrad((r\ball)^\circ)) \cdot |r\ball|.
     \end{eqnarray*} Similarly, for $\phi\in\Psi$, one has      \begin{eqnarray*}
n \phi (1/r) |r\ball| &\leq&  \OrliczG(r\ball)\leq \OrliczA(r\ball)  = n   \sup_{L\in \cS_0}V_{\phi}(r\ball , \vrad(L) L^\circ)   \\ &\leq &  n|r\ball| \sup_{L\in \cS_0}F\left( \frac{1}{n}\int _{S^{n-1}}\frac{r^n\r_L^n(u)}{|L|} \,d\s(u) \right)  =   n\phi\left( 1/r \right) |r\ball|,
\end{eqnarray*} where the last inequality follows from Jensen's inequality (for concave function $F$ as $\phi\in \Psi$). In conclusion, for  $\phi\in \Psi$, one also has  \begin{eqnarray*}
\OrliczG(r\ball)= \OrliczA(r\ball)=   n \phi\left(1/r \right) |r\ball|=n\phi (\vrad((r\ball)^\circ)) \cdot |r\ball|. \end{eqnarray*} 
Let $\E$ be any origin-symmetric ellipsoid with $|\E|=|r\ball|$ for some $r>0$. Then, $\E=T(r\ball)$ for some $T\in SL(n)$. Proposition \ref{homogeneous:degree} implies that, for $\phi\in \Phi$ or $\phi\in \Psi$,
\begin{eqnarray*} 
\OrliczA(\E) &=& \OrliczA(r\ball)=n \phi (\vrad((r\ball)^\circ)) \cdot |r\ball|=  n |\E|\cdot  \phi\left( \vrad(\E^\circ) \right); \\
  \OrliczG(\E)&=&\OrliczG(r\ball)=n\phi (\vrad((r\ball)^\circ)) \cdot |r\ball|= n  |\E| \cdot \phi\left(  \vrad(\E^\circ) \right).
      \end{eqnarray*} 
  In particular, for all $\E$ with $|\E|=|\ball|$, we have  \begin{eqnarray*}
 \OrliczG( \E)= \OrliczA( \E)=  n| \E| \cdot  \phi\left( 1 \right).
  \end{eqnarray*} 

The following proposition compares the Orlicz $L_{\phi}$ affine and geominimal surface areas with the Orlicz $\phi$-surface area of $K$ defined by $S_{\phi}(K)=nV_{\phi} (K, \ball)$. Clearly, \begin{eqnarray*} S_{\phi}(\ball)=nV_{\phi} (\ball, \ball)=\int _{S^{n-1}}\phi\left(1 \right) \,d\s(u) = \phi(1)\cdot n|\ball|= \OrliczG( \ball)= \OrliczA( \ball). \end{eqnarray*}  
\bp \label{Comparison:geominimal:affine:surface} Let $K\in \cK_0$ be a convex body in $\bbR^n$ with the origin in its interior. \\ (i) For   $\phi\in \Phi$, one has 
$$\OrliczA(K)\leq \OrliczG(K)\leq S_{\phi}(K).$$
(ii) For   $\phi\in \Psi$, one has 
$$\OrliczA(K)\geq \OrliczG(K)\geq S_{\phi}(K).$$\ep

\noindent {\bf Proof.} (i). Let $\phi\in \Phi$. Taking $Q=\ball$ in formula (\ref{Orlicz-geominimal-surface-area-convex}) and together with Proposition \ref{Comparison:geominimal:affine}, one has 
\begin{eqnarray*}
\OrliczA(K)\leq \OrliczG(K)= \inf_{Q\in \cK_0} \left\{nV_{\phi}(K, \vrad(Q^\circ)Q) \right\} \leq  nV_{\phi}(K, \ball) =S_{\phi}(K).\end{eqnarray*}
(ii). Let $\phi\in \Psi$. Taking $Q=\ball$ in formula (\ref{Orlicz-geominimal-surface-area:concave}) and together with Proposition \ref{Comparison:geominimal:affine}, one has 
\begin{eqnarray*}
\OrliczA(K) \geq \OrliczG(K)= \sup_{Q\in \cK_0} \left\{nV_{\phi}(K, \vrad(Q^\circ)Q) \right\}  \geq  nV_{\phi}(K, \ball) =S_{\phi}(K).\end{eqnarray*}
\bp \label{bounded by volume product} Let $K\in \cK_0$ be a convex body with the origin in its interior. 
\\ (i)  For $\phi\in \Phi$, one has $$\OrliczA(K)\leq \OrliczG(K) \leq \phi(\vrad(K^\circ)) \cdot n|K|.$$ 
(ii)  For $\phi\in \Psi$, one has $$\OrliczA(K)\geq \OrliczG(K) \geq \phi(\vrad(K^\circ)) \cdot n|K|.$$\ep
\noindent {\bf Proof.}  (i). Let $\phi\in \Phi$. Taking $Q=K$ in formula (\ref{Orlicz-geominimal-surface-area-convex}) and by Proposition \ref{Comparison:geominimal:affine}, one gets 
 $$\OrliczA(K)\leq \OrliczG(K) \leq nV_{\phi}(K,\vrad(K^\circ) K)=\phi(\vrad(K^\circ)) \cdot n|K|,$$ where for $K\in \cK_0$ and for  $\phi: (0, \infty)\rightarrow (0, \infty)$,   $$nV_{\phi}(K,\vrad(K^\circ)K)=\int_{S^{n-1}} \phi(\vrad(K^\circ)) h_K(u)\,dS(K, u)=\phi(\vrad(K^\circ))\cdot n|K|.$$
 
\noindent (ii).  Let $\phi\in \Psi$. Taking $Q=K$ in formula (\ref{Orlicz-geominimal-surface-area:concave}) and by Proposition \ref{Comparison:geominimal:affine}, one gets 
 $$\OrliczA(K)\geq \OrliczG(K) \geq nV_{\phi}(K,\vrad(K^\circ) K)=\phi(\vrad(K^\circ)) \cdot n|K|.$$ 
   \noindent {\bf Remark.} Replacing $K$ by its polar body $K^\circ$, one has, for $\phi\in \Phi$,
$$\OrliczA(K^\circ)\leq \OrliczG(K^\circ) \leq nV_{\phi}(K^\circ,\vrad(K) K^\circ)=\phi(\vrad(K)) \cdot n|K^\circ|.$$ Hence, for $\phi\in\Phi$, \begin{eqnarray*}
\OrliczA(K)\OrliczA(K^\circ) \leq  \OrliczG(K)\OrliczG(K^\circ)  \leq  \phi(\vrad(K))\phi(\vrad(K^\circ)) \cdot n^2|K|\cdot |K^\circ|.
\end{eqnarray*}  Moreover, if $\phi(t)\phi(s)\leq [\phi(1)]^2$ for all $t, s> 0$ satisfying $st\leq 1$, the following Santal\'{o} style inequality holds: for all $K\in \cK_s$ or $K\in \cK_c$,  \begin{eqnarray}
\OrliczA(K)\OrliczA(K^\circ)&\leq& \OrliczG(K)\OrliczG(K^\circ)\nonumber  \leq  \phi(\vrad(K))\phi(\vrad(K^\circ)) \cdot n^2|K|\cdot |K^\circ|\nonumber\\   &\leq & \phi^2(1) n^2\o_n^2 =[\OrliczG(\ball)]^2=[\OrliczA(\ball)]^2,\label{santalo style orlicz--1}
\end{eqnarray} where the last inequality follows from the Blaschke-Santal\'{o} inequality and the first equality follows from Corollary \ref{Orlicz:ellpsoids}. For instance, if $\phi(t)=t^p$ with $p\geq 0$, one gets the Santal\'{o} style inequality for $L_p$ affine and geominimal surface areas (see e.g.,  \cite{Lu1, WY2008, Ye2014, Zhou2011}). 

Similarly, for $\phi\in \Psi$ with $\phi(t)\phi(s)\geq A[\phi(1)]^2$ for some constant $A>0$ and for all $t, s> 0$ satisfying $st\leq 1$, one has, for all $K\in \cK_s$ or $K\in \cK_c$,  \begin{eqnarray}
\OrliczA(K)\OrliczA(K^\circ)&\geq& \OrliczG(K)\OrliczG(K^\circ)  \geq  \phi(\vrad(K))\phi(\vrad(K^\circ)) \cdot n^2|K|\cdot |K^\circ|\nonumber\\ &\geq& A c^n \phi^2(1) n^2\o_n^2 =Ac^n[\OrliczG(\ball)]^2=Ac^n[\OrliczA(\ball)]^2,\ \ \ \ 
\nonumber \end{eqnarray}  where the last inequality follows from the Bourgain-Milman inverse Santal\'{o} inequality and the first equality follows from Corollary \ref{Orlicz:ellpsoids}. For instance, if $\phi(t)=t^p$ with $-n<p<0$, one gets the Santal\'{o} style inequality for $L_p$ affine and geominimal surface areas \cite{WY2008, Ye2014}.

 The following theorem deals with the monotonicity for Orlicz $L_{\phi}$ affine and geominimal surface areas. We always assume that the function $\psi$ has inverse function $\psi^{-1}(\cdot)$. Let $H(t)=(\phi\circ \psi^{-1})(t)$ be the composition function of $\phi$ and $\psi^{-1}$. Recall that all functions  $\phi(t)\in \Psi$ are decreasing; and hence if $\phi(t), \psi(t)\in \Psi$, then $H(t)$ is increasing. Similar to the definition for $\Psi$, we are not interested in the case $H(t)$ being concave decreasing (as otherwise $H(t)$ is eventually a constant function, which leads to  $\phi$ being a constant function). Moreover, let $H(0)=\lim_{t\rightarrow 0} H(t)$ if the limit exists and is finite; on the other hand, let $H(0)=\infty$ if $\lim_{t\rightarrow 0} H(t)=\infty$. Similarly, one can also let $H(\infty)=\lim_{t\rightarrow \infty} H(t)$ if the limit exists and is finite; or simply $H(\infty)=\infty$ if $\lim_{t\rightarrow \infty} H(t)=\infty$.

\bt\label{cyclic}
Let $K\in \cK_0$ be a convex body with the origin in its interior. 
\vskip 2mm \noindent
(i) Assume that $\phi$ and $ \psi$ satisfy one of the following conditions: (a)  $\phi\in \Phi$ and $\psi\in \Psi$ with $H(t)$ increasing; (b) $\phi, \psi\in \Phi$ with $H(t)$ decreasing; (c) $H(t)$ concave increasing with either $\phi, \psi\in \Phi$ or $\phi, \psi\in \Psi$. Then \begin{equation*}
\frac{\OrliczA(K)}{n|K|} \leq  H\bigg(\frac{\O_{\psi}^{orlicz}(K)}{n|K|} \bigg)   \  \ and  \ \ \ \frac{\OrliczG(K)}{n|K|} \leq  H\bigg(\frac{G_{\psi}^{orlicz}(K)}{n|K|} \bigg) . \end{equation*}
(ii) Assume that $\phi$ and $ \psi$ satisfy one of the following conditions: (d) $\phi\in \Psi$ and $\psi\in \Phi$ with $H(t)$ increasing;   (e) $H(t)$ convex decreasing with one in $\Phi$ and another one in $\Psi$; (f) $H(t)$  convex increasing with either $\phi, \psi\in \Phi$ or $\phi, \psi\in \Psi$. Then \begin{equation*} 
\frac{\OrliczA(K)}{n|K|} \geq  H\bigg(\frac{\O_{\psi}^{orlicz}(K)}{n|K|} \bigg) \  \  and  \ \ \ \frac{\OrliczG(K)}{n|K|} \geq  H\bigg(\frac{G_{\psi}^{orlicz}(K)}{n|K|} \bigg) . \end{equation*}\et
  \noindent {\bf Remark.} When $\phi(t)=t^p$ and $\psi(t)=t^q$ for some $p, q \neq -n, 0$, one recovers the monotonicity properties for the $L_p$ affine and geominimal surface areas (see e.g., \cite{Lu1, WY2008, Ye2014, Zhang2007}). Clearly, condition (a) is equivalent to condition (d) if both $\phi(t)$ and $\psi(t)$ have inverse functions. If $H(t)$ is increasing and $\phi^{-1}(t), \psi^{-1}(t)$ both exist, condition (c) is equivalent to condition (f). This follows from the following claim: if $H(t)$ (and hence $H^{-1}(t)$) is increasing, $H(t)$ and $H^{-1}(t)$ have different convexity and concavity. In fact, without loss of generality, assume that $H(t)$ is increasing and convex. For all $t, s\in (0, \infty)$ and $\lambda\in [0, 1]$, we have, \begin{eqnarray*}
H^{-1}(\lambda t+(1-\lambda)s)&=&H^{-1}\big(\lambda H(H^{-1}(t))+(1-\lambda)H(H^{-1}(s))\big)\\ &\geq& H^{-1}\big( H(\lambda H^{-1}(t)+(1-\lambda) H^{-1}(s))\big)= \lambda H^{-1}(t)+(1-\lambda) H^{-1}(s),
\end{eqnarray*}  which leads to the concavity of $H^{-1}(t)$. Along the same line, one can also prove that the concavity of $H(t)$ implies the convexity of $H^{-1}(t)$. (A similar argument shows that if $H$ is decreasing, then $H(t)$ and $H^{-1}(t)$ have the same convexity and concavity.)

\vskip 2mm \noindent {\bf Proof.} We only prove the case for Orlicz $L_{\phi}$ geominimal surface area.  
 The proof for $\OrliczA(K)$ goes along the same line and hence is omitted.   
 
 \vskip 2mm \noindent (i). For condition (a) $\phi\in \Phi$ and $\psi\in \Psi$ with $H(t)$ increasing:  Proposition \ref{bounded by volume product} and the increasing property of $H(t)$ imply that 
\begin{eqnarray*} \frac{ \OrliczG(K)}{ n|K|}  \leq \phi(\vrad(K^\circ)) =H(\psi (\vrad(K^\circ))) \leq  H\bigg(\frac{ G_{\psi}^{orlicz}(K)}{ n|K|} \bigg). \end{eqnarray*}   
 For condition (b)  $\phi, \psi\in \Phi$ with $H(t)$ decreasing: Proposition \ref{bounded by volume product} and the decreasing property of $H(t)$ imply that 
\begin{eqnarray*} \frac{ \OrliczG(K)}{ n|K|}  \leq \phi(\vrad(K^\circ)) =H(\psi (\vrad(K^\circ))) \leq  H\bigg(\frac{ G_{\psi}^{orlicz}(K)}{ n|K|} \bigg). \end{eqnarray*}    
For condition (c): The concavity of $H(t)$ with Jensen's inequality imply that, $\forall K, Q\in \cK_0$,  \begin{eqnarray}\frac{V_{\phi}(K, \vrad(Q^\circ)Q)}{|K|}&=&\frac{1}{n|K|} \int _{S^{n-1}}H\bigg[\psi\left(\frac{\vrad(Q^\circ) h_Q(u)}{h_K(u)}\right)\bigg] h_K(u)\,dS(K, u)\nonumber\\ &\leq&H\bigg[ \frac{1}{n|K|} \int _{S^{n-1}}\psi\left(\frac{\vrad(Q^\circ) h_Q(u)}{h_K(u)}\right)h_K(u)\,dS(K, u) \bigg]\nonumber \\ &=& H\bigg(\frac{V_{\psi}(K, \vrad(Q^\circ)Q)}{|K|}\bigg)\nonumber.\end{eqnarray} Let $H(t)$ be increasing and concave with $\phi,\psi\in \Phi$. By formula (\ref{Orlicz-geominimal-surface-area-convex}), one has 
\begin{eqnarray}\frac{\OrliczG(K)}{n|K|} &=&\inf_{Q\in \cK_0} \frac{nV_{\phi}(K, \vrad(Q^\circ)Q)}{n|K|} \leq  \inf_{Q\in \cK_0}  H\bigg(\frac{nV_{\psi}(K, \vrad(Q^\circ)Q)}{n|K|}\bigg)\nonumber \\&=&   H\bigg(\inf_{Q\in \cK_0} \frac{nV_{\psi}(K, \vrad(Q^\circ)Q)}{n|K|}\bigg)=H\bigg(\frac{G_{\psi}^{orlicz}(K)}{n|K|} \bigg)\nonumber.\end{eqnarray}
Let $H(t)$ be increasing and concave with $\phi,\psi\in \Psi$. By formula (\ref{Orlicz-geominimal-surface-area:concave}), one has 
\begin{eqnarray}\frac{\OrliczG(K)}{n|K|} =\sup_{Q\in \cK_0} \frac{nV_{\phi}(K, \vrad(Q^\circ)Q)}{n|K|} \leq    H\bigg(\sup_{Q\in \cK_0} \frac{nV_{\psi}(K, \vrad(Q^\circ)Q)}{n|K|}\bigg)=H\bigg(\frac{G_{\psi}^{orlicz}(K)}{n|K|} \bigg)\nonumber.\end{eqnarray}
 
 \noindent (ii). For condition (d) $\phi\in \Psi$ and $\psi\in \Phi$ with $H(t)$ increasing: Proposition \ref{bounded by volume product} and the increasing property of $H(t)$ imply that 
\begin{eqnarray*} \frac{ \OrliczG(K)}{ n|K|}  \geq \phi(\vrad(K^\circ)) =H(\psi (\vrad(K^\circ))) \geq  H\bigg(\frac{ G_{\psi}^{orlicz}(K)}{ n|K|} \bigg).\end{eqnarray*}  
For condition (e): the convexity of $H(t)$ together with Jensen's inequality imply that \begin{eqnarray}\frac{V_{\phi}(K, \vrad(Q^\circ)Q)}{|K|}&=& \frac{1}{n|K|} \int _{S^{n-1}}H\bigg[\psi\left(\frac{\vrad(Q^\circ) h_Q(u)}{h_K(u)}\right)\bigg] h_K(u)\,dS(K, u)\nonumber\\ &\geq& H\bigg(\frac{V_{\psi}(K, \vrad(Q^\circ)Q)}{|K|}\bigg) \label{convex:decreasing---1}.\end{eqnarray}
Let $\phi \in \Psi$ and $\psi\in \Phi$ with $H(t)$ convex decreasing: taking the supremum over $Q\in \cK_0$, together with the decreasing property of $H(t)$ and formulas (\ref{Orlicz-geominimal-surface-area-convex}) and (\ref{Orlicz-geominimal-surface-area:concave}), one has 
 \begin{eqnarray*}\frac{\OrliczG(K)}{n|K|}   \geq   \sup_{Q\in \cK_0}  H\bigg(\frac{nV_{\psi}(K, \vrad(Q^\circ)Q)}{n|K|}\bigg) =   H
 \bigg(\inf_{Q\in \cK_0} \frac{nV_{\psi}(K, \vrad(Q^\circ)Q)}{n|K|}\bigg)=H\bigg(\frac{G_{\psi}^{orlicz}(K)}{n|K|} \bigg).\end{eqnarray*} 
 Similarly, for $\phi \in \Phi$ and $\psi\in \Psi$ with $H(t)$ convex decreasing: 
  \begin{eqnarray*}\frac{\OrliczG(K)}{n|K|}   \geq   \inf_{Q\in \cK_0}  H\bigg(\frac{nV_{\psi}(K, \vrad(Q^\circ)Q)}{n|K|}\bigg) =   H
  \bigg(\sup_{Q\in \cK_0} \frac{nV_{\psi}(K, \vrad(Q^\circ)Q)}{n|K|}\bigg)=H\bigg(\frac{G_{\psi}^{orlicz}(K)}{n|K|} \bigg).\end{eqnarray*} 
For condition (f) $\phi, \psi\in \Phi$ with $H(t)$ convex increasing: taking the infimum over $Q\in \cK_0$ in inequality (\ref{convex:decreasing---1}), together with formula (\ref{Orlicz-geominimal-surface-area-convex}), one has 
\begin{eqnarray*}\frac{\OrliczG(K)}{n|K|}   \geq   \inf_{Q\in \cK_0}  H\bigg(\frac{nV_{\psi}(K, \vrad(Q^\circ)Q)}{n|K|}\bigg) =   H
 \bigg(\inf_{Q\in \cK_0} \frac{nV_{\psi}(K, \vrad(Q^\circ)Q)}{n|K|}\bigg)=H\bigg(\frac{G_{\psi}^{orlicz}(K)}{n|K|} \bigg).\end{eqnarray*} 
 Similarly, for $\phi, \psi\in \Psi$ with $H(t)$ convex increasing: taking the supremum over $Q\in \cK_0$ in inequality (\ref{convex:decreasing---1}), together with formula (\ref{Orlicz-geominimal-surface-area:concave}), one has 
\begin{eqnarray*}\frac{\OrliczG(K)}{n|K|}   \geq   \sup_{Q\in \cK_0}  H\bigg(\frac{nV_{\psi}(K, \vrad(Q^\circ)Q)}{n|K|}\bigg) =   H
 \bigg(\sup_{Q\in \cK_0} \frac{nV_{\psi}(K, \vrad(Q^\circ)Q)}{n|K|}\bigg)=H\bigg(\frac{G_{\psi}^{orlicz}(K)}{n|K|} \bigg).\end{eqnarray*} 

Let $\psi=t\in \Phi$ and $H(t)=\phi(t)\in \Phi$ be concave increasing. Part (i) of Theorem \ref{cyclic} implies \begin{equation}\label{i-2-removing-1}  
\frac{\OrliczA(K)}{n|K|} \leq  \phi\bigg(\frac{\O_{1}^{orlicz}(K)}{n|K|} \bigg) \ \ \ and  \ \ \ \frac{\OrliczG(K)}{n|K|} \leq  \phi \bigg(\frac{G_{1}^{orlicz}(K)}{n|K|} \bigg) . \end{equation}
If $\psi=t\in \Phi$ and $\phi(t)=H(t)\in \Psi$ is convex decreasing, part (ii) of Theorem \ref{cyclic}  implies  
\begin{equation} 
\frac{\OrliczA(K)}{n|K|} \geq  \phi \bigg(\frac{\O_{1}^{orlicz}(K)}{n|K|} \bigg) \ \ \  and  \ \ \ \frac{\OrliczG (K)}{n|K|} \geq  \phi\bigg(\frac{G_{1}^{orlicz}(K)}{n|K|} \bigg).  \label{i-2-removing-2}  \end{equation}  
  
  We now prove the following affine isoperimetric inequality for Orlicz $L_{\phi}$ affine and geominimal surface areas. Write $B_K$ for the origin-symmetric Euclidean ball with $|B_K|=|K|$ for $K\in \cK_0$. 

\bt \label{isoperimetric:geominimal}
Let $K\in \cK_c$ or $K\in \cK_s$. 

\vskip 2mm \noindent (i) For $\phi\in \Phi$, the following affine isoperimetric inequality holds, with equality if and only if $K$ is an origin-symmetric ellipsoid, 
\begin{eqnarray*}
\OrliczA(K)\leq \OrliczG(K)\leq \OrliczG([B_{K^\circ}]^\circ)= \OrliczA([B_{K^\circ}]^\circ). 
\end{eqnarray*} 
If in addition $\phi$ is increasing, then
\begin{eqnarray*}
\OrliczA(K)\leq \OrliczG(K)\leq \OrliczG(B_K)= \OrliczA(B_K),
\end{eqnarray*} with equality if $K$ is an origin-symmetric ellipsoid. Moreover, if $\phi$ is strictly increasing, equality holds if and only if $K$ is an origin-symmetric ellipsoid. 

\vskip 2mm \noindent 
(ii) For $\phi\in \Psi$, the following affine isoperimetric inequality holds with equality if $K$ is an origin-symmetric ellipsoid,
\begin{eqnarray*}
\OrliczA(K)\geq \OrliczG(K)\geq \OrliczG(B_K)= \OrliczA(B_K).
\end{eqnarray*}  If in addition $\phi$ is strictly decreasing, equality holds if and only if $K$ is an origin-symmetric ellipsoid.  \et 

 \noindent {\bf Proof.} First, Blaschke-Santal\'{o} inequality implies that, for all $K\in \cK_c$ or $K\in \cK_s$, $$\vrad(K)\vrad(K^\circ)\leq 1=\vrad(B_K)\vrad([B_K]^\circ)=\vrad(K)\vrad([B_K]^\circ).$$ Therefore, for all  $K\in \cK_c$ or $K\in \cK_s$, one has \begin{equation}
  \label{compare:volume:radius:K} \vrad(K^\circ)\leq \vrad([B_K]^\circ), 
  \end{equation} with equality if and only if $K$ is an origin-symmetric ellipsoid.

 \vskip 2mm \noindent  (i). Let $\phi\in \Phi$. Corollary \ref{Orlicz:ellpsoids} and Proposition  \ref{bounded by volume product} imply 
   \begin{eqnarray*}
  \OrliczA(K)&\leq& \OrliczG(K) \leq \phi(\vrad(K^\circ)) \cdot n|K|= \phi(\vrad(B_{K^\circ})) \cdot n|[B_{K^\circ}]^\circ|\cdot \frac{|K|}{|[B_{K^\circ}]^\circ| } \\ &=& 
  \OrliczG([B_{K^\circ}]^\circ)\frac{|K||K^\circ|}{|[B_{K^\circ}]^\circ| |B_{K^\circ}| }  \leq   
 \OrliczG([B_{K^\circ}]^\circ)=\OrliczA([B_{K^\circ}]^\circ), 
  \end{eqnarray*} where the last inequality follows from the Blaschke-Santal\'{o} inequality. 
  Clearly, equality holds if $K$ is an origin-symmetric ellipsoid. On the other hand, equality holds only if equality holds in the Blaschke-Santal\'{o} inequality, that is, $K$ has to be an origin-symmetric ellipsoid.

 Assume in addition that $\phi\in \Phi$ is increasing. Inequality (\ref{compare:volume:radius:K}) together with Corollary \ref{Orlicz:ellpsoids} and Proposition \ref{bounded by volume product} imply 
  \begin{eqnarray*}
 \OrliczA(K) \leq  \OrliczG(K) \leq n|K|  \phi(\vrad(K^\circ))  \leq   n|B_K|  \phi(\vrad([B_K]^\circ))=\OrliczG(B_K)= \OrliczA(B_K),\end{eqnarray*} 
 with equality if $K$ is an origin-symmetric ellipsoid. Moreover, if $\phi$ is strictly increasing, equality holds only if equality holds in inequality (\ref{compare:volume:radius:K}), i.e., $K$ has to be an origin-symmetric ellipsoid.

\vskip 2mm \noindent  (ii). Let $\phi\in \Psi$. Recall that $\phi\in \Psi$ is decreasing. Inequality  (\ref{compare:volume:radius:K}) together with Corollary \ref{Orlicz:ellpsoids} and Proposition \ref{bounded by volume product} imply 
  \begin{eqnarray*}
 \OrliczA(K) \geq  \OrliczG(K) \geq n|K| \phi(\vrad(K^\circ))  \geq   n|B_K|  \phi(\vrad([B_K]^\circ) ) =\OrliczG(B_K)= \OrliczA(B_K),  
 \end{eqnarray*} 
 with equality if $K$ is an origin-symmetric ellipsoid. Moreover if $\phi$ is strictly decreasing, equality holds only if equality holds in inequality (\ref{compare:volume:radius:K}), i.e., $K$ has to be an origin-symmetric ellipsoid.

\vskip 2mm \noindent {\bf Remark.}  Part (i) of Theorem \ref{isoperimetric:geominimal} asserts that {\em among all convex bodies in $\cK_c$ (or $\cK_s$) with fixed volume, Orlicz $L_{\phi}$ affine and geominimal surface areas for $\phi\in \Phi$ attain the maximum at origin-symmetric ellipsoids.} Moreover, if $\phi$ is also strictly increasing, the origin-symmetric ellipsoids are the only maximizers. 
  Similarly, part (ii) of Theorem \ref{isoperimetric:geominimal} asserts that {\em among all convex bodies in $\cK_c$ (or $\cK_s$) with fixed volume,   Orlicz $L_{\phi}$ affine and geominimal surface areas for strictly decreasing function $\phi\in \Psi$ attain the minimum at and only at origin-symmetric ellipsoids.} When $\phi(t)=t^p$ with $-n\neq p\in \bbR$, one recovers the $L_p$ affine isoperimetric inequalities for the $L_p$ affine and geominimal surface areas (see e.g., \cite{Lu1, Petty1974, Petty1985, WY2008, Ye2014}). 

The following result removes the condition on the centroid (or the Santal\'{o} point) of $K$ in Theorem \ref{isoperimetric:geominimal} for certain $\phi$.  The  case $\phi(t)=t^p$ with $p\in (-n, 1)$ were proved in \cite{Ye2014, Zhang2007}. See also \cite{Ye2013} for similar results on general affine surface areas. 
 \bc \label{Lp:p:positive:symmetrization} Let $K\in\cK_0$ be a convex body with the origin in its interior. 
\\
\noindent (i). Let $\phi\in \Phi$ be a concave increasing function. The Orlicz $L_{\phi}$ affine and geominimal surface areas both attain the maximum at origin-symmetric ellipsoids among all convex bodies in $\cK_0$ with fixed volume. More
precisely, \begin{eqnarray*} \OrliczA(K)\leq\OrliczG(K)\leq \OrliczG(B_K)= \OrliczA(B_K)\end{eqnarray*} with equality if $K$ is an origin-symmetric ellipsoid. Moreover, if $\phi$ is increasing and strictly concave,  equality holds if and only if $K$ is an origin-symmetric ellipsoid. \\ \noindent 
(ii). Let $\phi\in \Psi$ be a convex decreasing function. The Orlicz $L_{\phi}$ affine and geominimal surface areas both attain the minimum at origin-symmetric ellipsoids among all convex bodies in $\cK_0$ with fixed volume. More
precisely, \begin{eqnarray*} \OrliczA(K)\geq\OrliczG(K)\geq \OrliczG(B_K)= \OrliczA(B_K)\end{eqnarray*} with equality if $K$ is an origin-symmetric ellipsoid. Moreover, if $\phi$ is decreasing and strictly convex, equality holds if and only if $K$ is an origin-symmetric ellipsoid.  \ec
  \noindent {\bf Proof.} We only prove the case for Orlicz $L_{\phi}$ geominimal surface area. The proof for Orlicz $L_{\phi}$ affine surface area goes along the same line.

\vskip 2mm \noindent  (i). Let $\phi\in \Phi$ be a concave increasing function. Assume that $z_0$ is the centroid of $K$. By inequality   (\ref{i-2-removing-1})  and the translation invariance of $ G_{1}^{orlicz}(K)$, one has,  \begin{eqnarray}
 \frac{\OrliczG(K)}{n|K|}  \!\leq\!   \phi \bigg(\!\frac{G_{1}^{orlicz}(K)}{n|K|}\! \bigg) \!=\! \phi \bigg(\!\frac{G_{1}^{orlicz}(K-z_0)}{n|K-z_0|}\! \bigg)\! \leq \! \phi \bigg(\!\frac{G_{1}^{orlicz}(B_K)}{n|B_K|}\! \bigg)\!=\!\phi(\vrad([B_K]^\circ)), \ \ \ \label{isoperimetric:remove-1} \end{eqnarray} where the second inequality follows from Theorem \ref{isoperimetric:geominimal} (as $K-z_0\in \cK_c$) and the last equality follows from Corollary \ref{Orlicz:ellpsoids} (with function $\phi(t)=t$). Therefore, again by Corollary \ref{Orlicz:ellpsoids}, one has  \begin{eqnarray*}
  {\OrliczG(K)}\leq  \phi(\vrad([B_K]^\circ))\cdot  {n|K|} =\phi(\vrad([B_K]^\circ)) \cdot {n|B_K|} ={\OrliczG(B_K)}. \end{eqnarray*} Clearly, equality holds if $K$ is an origin-symmetric ellipsoid. 
  
  Assume in addition that $\phi$ is strictly concave. To have equality in inequality (\ref{isoperimetric:remove-1}), one requires $K-z_0$ to be an origin-symmetric ellipsoid. We now claim that $z_0=0$. To this end, assume that, without loss of generality,  $K=z_0+r\ball$ is a Euclidean ball with center $z_0$ and radius $r>0$. Moreover, $\|z_0\|<r$ and $$\int _{S^{n-1}}  {h_{z_0+r\ball}(u)}  \,d\s(u) =nr|\ball|.$$ Hence, for strictly concave function $\phi\in\Phi$,  \begin{eqnarray}\OrliczG(K) &=& \OrliczG(z_0+r\ball )\leq nV_{\phi}(z_0+r\ball, \ball)\nonumber \\  &=& n\o_n r^n\int _{S^{n-1}} \phi\bigg(\!\frac{1}{h_{z_0+r\ball}(u)}\!\bigg) \frac{h_{z_0+r\ball}(u)}{nr\o_n} \,d\s(u)\nonumber\\  &\leq& n\o_n r^n \cdot  \phi\bigg(\int _{S^{n-1}} \frac{h_{z_0+r\ball}(u)}{nr\o_n h_{z_0+r\ball}(u)} \,d\s(u) \bigg)\nonumber \\ &=& n \o_nr^n \cdot \phi(1/r)=\OrliczG(r\ball), \label{non translation invariance 1} \end{eqnarray} where the last inequality follows from Jensen's inequality (for concave function $\phi$) and the last equality follows from Corollary \ref{Orlicz:ellpsoids}. To have equality in inequality (\ref{non translation invariance 1}), one requires, in particular, equality for Jensen's  inequality. That is, $h_{z_0+r\ball}(u)$ is a constant on $S^{n-1}$ due to the strict concavity of $\phi$, and hence $z_0=0$ as desired.

\vskip 2mm \noindent (ii). Let  $\phi\in \Psi$ be convex decreasing. Assume that $z_0$ is the centroid of $K$. By inequality   (\ref{i-2-removing-2})  and the translation invariance of $ G_{1}^{orlicz}(K)$, one has,  \begin{eqnarray*} \phi^{-1}\bigg(\frac{\OrliczG(K)}{n|K|} \bigg) \leq   \frac{G_{1}^{orlicz}(K)}{n|K|}= \frac{G_{1}^{orlicz}(K-z_0)}{n|K-z_0|}  \leq   \frac{G_{1}^{orlicz}(B_K)}{n|B_K|}  = \vrad([B_K]^\circ) , \label{isoperimetric:remove-2} \end{eqnarray*} where the second inequality follows from Theorem \ref{isoperimetric:geominimal} (as $K-z_0\in \cK_c$) and the last equality follows from Corollary \ref{Orlicz:ellpsoids} (with function $\phi(t)=t$). As $\phi$ is decreasing and by Corollary \ref{Orlicz:ellpsoids}, one has  \begin{eqnarray}\frac{\OrliczG(K)}{n|K|}   \geq   \phi \big( \vrad([B_K]^\circ)\big) \Longleftrightarrow {\OrliczG(K)}   \geq  \phi \big( \vrad([B_K]^\circ)\big) {n|K|}=\OrliczG(B_K). \label{remove-1-1-1} \end{eqnarray} Clearly, equality holds if $K$ is an origin-symmetric ellipsoid.

  Assume in addition that $\phi$ is convex and strictly decreasing. To have equality in inequality (\ref{remove-1-1-1}), one requires $K-z_0$ to be an origin-symmetric ellipsoid, which implies $z_0=0$. In fact, similar to the case (i), one has \begin{eqnarray}\OrliczG(K) &=& \OrliczG(z_0+r\ball )=\sup_{Q\in \cK_0} \left\{nV_{\phi}(z_0+r\ball, \vrad(Q^\circ)Q) \right\} \nonumber \\  &\geq&  n\o_n r^{n-1}\int _{S^{n-1}} \phi\bigg(\frac{1}{h_{z_0+r\ball}(u)}\bigg) \frac{h_{z_0+r\ball}(u)}{n\o_n} \,d\s(u)\nonumber \\ &\geq & n\o_nr^n \cdot \phi(1/r)=\OrliczG(r\ball). \ \ \ \ \ \ \label{non translation invariance 2} \end{eqnarray}  To have equality in (\ref{non translation invariance 2}), one requires, in particular, equality for Jensen's  inequality. That is, $h_{z_0+r\ball}(u)$ is a constant on $S^{n-1}$ due to the strict convexity of $\phi$, and hence $z_0=0$ as desired.


\section{Orlicz mixed $L_{\phi}$ affine and geominimal surface areas for multiple convex bodies} \label{section 5}
In this section, the Orlicz mixed $L_{\phi}$ affine and geominimal surface areas for multiple convex bodies and their basic properties are  briefly discussed.  We will omit most of the proofs because these proofs are either similar to those for single convex body discussed in Section \ref{section 3} or similar to those for mixed $p$-affine and mixed $L_p$ geominimal surface areas in \cite{Lu1, WernerYe2010, DYzhuzhou2014}.  

We use $\vec{\phi}$ for $(\phi_1, \phi_2, \cdots, \phi_n)$. We say $\vec{\phi}\in \Phi^n$ (or $\vec{\phi}\in \Psi^n$) if each $\phi_i\in \Phi$ (or $\phi_i\in \Psi$). Similarly, $\cL=(L_1, \cdots, L_n)\in \cS_0^n$ and $\bK=(K_1, \cdots, K_n)\in \cF_0^n$ mean that each $L_i\in \cS_0$ and each $K_i\in \cF_0$ respectively. We use  $\cL^\circ$ for $(L_1^\circ, \cdots, L_n^\circ)$. Define  $V_{\vec{\phi}}(\bK, \bQ)$ for $\bK\in \cF_0^n$ and $\bQ\in \cK_0^n$ by \be\nonumber V_{\vec{\phi}}(\bK, \bQ)=\frac{1}{n}\int _{S^{n-1}}\prod_{i=1}^n \bigg[\phi_i\left(\frac{h_{Q_i}(u)}{h_{K_i}(u)}\right) h_{K_i}(u)f_{K_i}(u)\bigg]^{\frac{1}{n}} \,d\s(u).\ee When $\bK\in \cF_0^n$ and $\cL\in \cS_0^n$, we use the following formula  $$V_{\vec{\phi}}(\bK, \cL^\circ)=\frac{1}{n}\int _{S^{n-1}}\prod_{i=1}^n \bigg[\phi_i\left(\frac{1}{\r_{L_i}(u)h_{K_i}(u)}\right) h_{K_i}(u)f_{K_i}(u)\bigg]^{\frac{1}{n}} \,d\s(u).$$ When $\phi_i=\phi$, 
$K_i=K$ and $L_i=L$ for all $i=1, 2, \cdots, n$, one gets $ V_{\vec{\phi}} (\bK; \cL^\circ)=V_{\phi}(K, L^\circ). $  

We now propose our definition for the Orlicz mixed $L_{\phi}$ affine surface area. 
\bd\label{equivalent:mixed affine:surface:area-1} Let $K_1, \cdots,
K_n\in \cF_0$. \\
(i) For $\vec{\phi}\in \Phi^n$, we define  $\OrliczAmix(\bK)$ by
 \begin{equation*} \OrliczAmix(\bK)= \inf_{ \cL\in \cS_0^n }
 \big\{n V_{\vec{\phi}}(\bK;  \cL^\circ) \ \ with \ \  |L_1|=|L_2|=\cdots =|L_n|=\o_n\big\}. \end{equation*}
(ii) For $\vec{\phi}\in \Psi^n$, we define  $\OrliczAmix(\bK)$ by
 \begin{equation*} \OrliczAmix(\bK)= \sup_{ \cL\in \cS_0^n }
 \big\{n V_{\vec{\phi}}(\bK;  \cL^\circ) \ \ with \ \  |L_1|=|L_2|=\cdots =|L_n|=\o_n\big\}. \end{equation*}
\ed
Similarly,  the Orlicz mixed $L_{\phi}$ geominimal surface area can be defined as follows. 
\bd\label{equivalent:mixed geominimal:surface:area-1}
Let $K_1, \cdots,
K_n\in \cF_0$. \\
(i) For $\vec{\phi}\in \Phi^n$, we define  $\OrliczGmix(\bK)$ by
 \begin{equation*} \OrliczGmix(\bK)= \inf_{ \bQ\in \cK_0^n }
 \big\{n V_{\vec{\phi}}(\bK;  \bQ) \ \ with \ \  |Q_1^\circ|=|Q_2^\circ|=\cdots =|Q_n^\circ|=\o_n\big\}. \end{equation*}
(ii) For $\vec{\phi}\in \Psi^n$, we define  $\OrliczGmix(\bK)$ by
 \begin{equation*} \OrliczGmix(\bK)= \sup_{ \bQ\in \cK_0^n }
 \big\{n V_{\vec{\phi}}(\bK;  \bQ) \ \ with \ \  |Q_1^\circ|=|Q_2^\circ|=\cdots =|Q_n^\circ|=\o_n\big\}. \end{equation*}\ed
 \noindent {\bf Remark.} The case $\vec{\phi}=(t^p, t^p, \cdots, t^p)$ corresponds to  the mixed $p$-affine and mixed $L_p$ geominimal surface areas (see e.g., \cite{Lu1, WernerYe2010, DYzhuzhou2014}). For the geominimal case, several different mixed $L_p$ geominimal surface areas could be proposed for the same $\bK$. Analogously, one can also define several mixed $L_{\phi}$ Orlicz affine and geominimal surface areas; however, due to high similarity of their properties (as one can see in \cite{DYzhuzhou2014}), we only focus on the one by Definition  \ref{equivalent:mixed geominimal:surface:area-1}  in this paper.

The following result could be proved by a similar argument to Proposition \ref{Comparison:geominimal:affine}.
 \bp\label{proposition 4.3--1}   Let $K_1, \cdots, K_n\in \cF_0$. \\ (i) For $\vec{\phi}\in \Phi^n$, one has 
$$\OrliczAmix(\bK)\leq \OrliczGmix(\bK).$$
(ii) For $\vec{\phi}\in \Psi^n$, one has 
$$\OrliczAmix(\bK)\geq \OrliczGmix(\bK).$$
 \ep

 The following result states that Orlicz mixed $L_{\phi}$ affine and geominimal surface areas are affine invariant. The proof is similar to that of Proposition \ref{homogeneous:degree}.  For $T\in SL(n)$ and $\bK=(K_1, \cdots, K_n)\in \cK_0^n$, we let $T \bK=(T K_1, \cdots, T K_n)$.  
 
 \bp \label{proposition-1} Let $\bK\in \cF_0^n$. For $\vec{\phi}\in \Phi^n$ or $\vec{\phi}\in \Psi^n$, one has $$ \OrliczAmix(T\bK)=\OrliczAmix(\bK); \ \ \ \OrliczGmix(T\bK)=\OrliczGmix(\bK), \ \  \forall T\in SL(n).$$ \ep 

 The classical Alexander-Fenchel inequality for the mixed volume  (see
 \cite{Sch}) is one of the most important inequalities in convex geometry. It has been extended to the mixed
$p$-affine surface area and the mixed $L_p$ geominimal surface area (see e.g., \cite{Lut1987, Lu1, WernerYe2010, DYzhuzhou2014, zhuzhou}). See also \cite{Y} for similar inequalities related to general mixed affine surface areas. Here, we prove the following Alexander-Fenchel type
inequality for Orlicz mixed $L_{\phi}$ affine and geominimal surface areas.

\bt \label{Alexander-Fenchel:both} Let $\bK\in \cF_0^n $. For $\vec{\phi}\in \Phi^n$ or $\vec{\phi}\in \Psi^n$, one has \begin{eqnarray*}  
\big[\OrliczAmix(\bK)\big]^{n} \leq \prod_{i=1}^n \O_{\phi_i}^{orlicz}(K_i) \ \ and \ \ \  \big[\OrliczGmix(\bK)\big]^{n}  \leq \prod_{i=1}^n  G_{\phi_i}^{orlicz}(K_i). \end{eqnarray*} Moreover, if $\vec{\phi}\in \Psi^n$, the following Alexander-Fenchel type inequality holds:  for $1 \leq m \leq n$,  
\begin{eqnarray*}
\big[\OrliczAmix(\bK)\big]^{m}&\leq& \prod_{i=0}^{m-1}\O_{(\phi_{1},
\cdots, \phi_{n-m}, \phi_{n-i},\cdots, \phi_{n-i})}^{orlicz}(K_{1},
\cdots, K_{n-m}, \underbrace{K_{n-i},\cdots, K_{n-i}}_{m}),  \\  \big[\OrliczGmix(\bK)\big]^{m}&\leq& \prod_{i=0}^{m-1}G_{(\phi_{1},
\cdots, \phi_{n-m}, \phi_{n-i},\cdots, \phi_{n-i})}^{orlicz}(K_{1},
\cdots, K_{n-m}, \underbrace{K_{n-i},\cdots, K_{n-i}}_{m}).
\end{eqnarray*}  \et
  \noindent {\bf Remark.} If each $\phi_i\in \Phi$ satisfies $\phi_i(t)\phi_i(s)\leq [\phi(1)]^2$ for all $t, s> 0$ with $st\leq 1$, then   \begin{eqnarray*} \big[\OrliczAmix(\bK)\OrliczAmix(\bK^\circ)\big]^n&\leq& \big[\OrliczGmix(\bK) \OrliczGmix(\bK^\circ)\big]^n  \leq  \prod_{i=1}^n\big[G_{\phi_i}^{orlicz}(K_i)G_{\phi_i}^{orlicz}(K_i^\circ)\big] \\    &\leq &  \prod_{i=1}^n [G_{\phi_i}^{orlicz}(\ball)]^2= \prod_{i=1}^n  [\O_{\phi_i}^{orlicz} (\ball)]^2,\label{santalo style orlicz-1}
\end{eqnarray*} where the last inequality follows from inequality (\ref{santalo style orlicz--1}).

 \noindent \textbf{Proof.} We only prove the geominimal case and omit the proof for affine case.  
Let \begin{eqnarray*} H_{0}(u)&=&\prod_{i=1}^{n-m} \bigg[\phi_i\left(\frac{h_{Q_i}(u)}{h_{K_i}(u)}\right) h_{K_i}(u)f_{K_i}(u)\bigg]^{\frac{1}{n}},\\ H_{i+1}(u)&=& \bigg[\phi_{n-i}\left(\frac{h_{Q_{n-i}}(u)}{h_{K_{n-i}}(u)}\right) h_{K_{n-i}}(u)f_{K_{n-i}}(u)\bigg]^{\frac{1}{n}}, \ \ i=0, \cdots, m-1. \end{eqnarray*}
H\"{o}lder's inequality (see \cite{HLP}) implies 
\begin{eqnarray}
 \big[V_{\vec{\phi}}(\bK;  \bQ)\big]^m 
 &=& \left(\frac{1}{n}\int_{S^{n-1}}H_{0}(u)H_{1}(u)\cdots H_{m}(u) d\s(u)\right)^m \nonumber \\
  &\leq& \prod_{i=0}^{m-1}\Big(\frac{1}{n}\int_{S^{n-1}}H_{0}(u)[H_{i+1}(u)]^md\s(u)\Big)  =\prod_{i=0}^{m-1} A_i,\label{Holder:mixed:-001} \end{eqnarray} where for $i=0, 1, \cdots, m-1$,  \begin{eqnarray*}A_i
=  V_{(\phi_{1}, \cdots , \phi_{n-m}, \phi_{n-i},\cdots, \phi_{n-i}) } (K_{1},\!\cdots\!,\! K_{n-m}, \underbrace{K_{n-i},\!\cdots\!, K_{n-i}}_{m}; Q_{1},\!\cdots\!, Q_{n-m}, \underbrace{Q_{n-i},\!\cdots\!, Q_{n-i}}_{m}). 
\end{eqnarray*} 
In particular, if $m=n$, inequality (\ref{Holder:mixed:-001}) implies that \begin{eqnarray} \nonumber
 \big[V_{\vec{\phi}}(\bK;  \bQ)\big]^n&\leq &\prod_{i=1}^{n} V_{\phi_i} (K_{i},  Q_{i}). 
\end{eqnarray} 
Let $\vec{\phi}\in \Phi^n$. By Definitions  \ref{Orlicz geominimal surface} and \ref{equivalent:mixed
 geominimal:surface:area-1},  one has 
\begin{eqnarray*}
\big[\OrliczGmix(\bK)\big]^n &=& \big[\inf_{ \bQ\in \cK_0^n }
 \big\{n V_{\vec{\phi}}(\bK;  \bQ) \ \ with \ \  |Q_1^\circ|=|Q_2^\circ|=\cdots =|Q_n^\circ|=\o_n\big\}\big]^n \\  &\leq &\inf_{\bQ\in \cK_0^n}  \bigg\{\prod_{i=1}^{n} [nV_{\phi_i} (K_{i},  Q_{i})] \ \ with \ \  |Q_1^\circ| =\cdots =|Q_n^\circ|=\o_n\bigg\} \\ &=&\prod_{i=1}^{n} \inf_{ Q_i\in \cK_0^n}  \big\{nV_{\phi_i} (K_{i},  Q_{i})\ \ with \ \  |Q_i^\circ| =\o_n\big\}=\prod_{i=1}^n G_{\phi_i}^{orlicz}(K_i).\end{eqnarray*}  
 Similarly, if $\phi\in \Psi^n$, one gets \begin{eqnarray*}
\big[\OrliczGmix(\bK)\big]^n  &\leq & \prod_{i=1}^{n} \sup_{ Q_i\in \cK_0^n}  \big\{nV_{\phi_i} (K_{i},  Q_{i})\ \ with \ \  |Q_i^\circ| =\o_n\big\}=\prod_{i=1}^n G_{\phi_i}^{orlicz}(K_i). \end{eqnarray*}  Moreover, inequality (\ref{Holder:mixed:-001}) implies that for $\vec{\phi}\in \Psi^n$ and $1\leq m\leq n$,
\begin{eqnarray*}  \big[\OrliczGmix(\bK)\big]^{m} 
 &\leq&
 \prod_{i=0}^{m-1}\sup_{\bQ \in \cK_0^n} \{ A_i\ \ with\ \  |Q_1^\circ|=\cdots =|Q_n^{\circ}|=\o_n\},\end{eqnarray*} and the desired Alexander-Fenchel type inequality follows if one notices $$\sup_{\bQ \in \cK_0^n} \{ A_i:   |Q_1^\circ|\!=\!\cdots\! =\! |Q_n^\circ|\!=\!\o_n\}\leq G_{(\phi_{1},
\cdots, \phi_{n-m}, \phi_{n-i},\cdots, \phi_{n-i})}^{orlicz}(K_{1},\!
\cdots\!, K_{n-m}, \underbrace{K_{n-i},\!\cdots\!, K_{n-i}}_{m}).$$

The following result is a direct consequence from Theorem \ref{Alexander-Fenchel:both} and Propositions  \ref{Comparison:geominimal:affine:surface} and \ref{proposition 4.3--1}. 
 \bc Let $K_1, \cdots, K_n\in \cF_0$.  For $\vec{\phi}\in \Phi^n$, one has \begin{eqnarray*}  
\big[\OrliczAmix(\bK)\big]^{n} \leq  \big[\OrliczGmix(\bK)\big]^{n} \leq  S_{\phi_1}(K_1)\cdots S_{\phi_n}(K_n). \end{eqnarray*}   \ec
 
 A direct consequence of Proposition \ref{proposition 4.3--1}, and  Theorems \ref{isoperimetric:geominimal} and \ref{Alexander-Fenchel:both} is the following isoperimetric type inequality, which holds for any possible combinations of $\cF_c$ and $\cF_s$. Due to Corollary \ref{Lp:p:positive:symmetrization}, $K_i$ can be assumed even in $\cF_0$ if $\phi_i \in \Phi$ is concave increasing.  

\bt \label{isoperimetric:geominimal:mixed case}
 Let $K_{i} \in \cF_c$ or $K_i\in \cF_s$. For $\vec{\phi}\in \Phi^n$, one has 
\begin{eqnarray*}
\big[\OrliczAmix(\bK)\big]^{n} \leq  \big[\OrliczGmix(\bK)\big]^{n} \leq \prod_{i=1}^n G^{orlicz}_{\phi_i}([B_{(K_i)^\circ}]^\circ)=\prod_{i=1}^n \O^{orlicz}_{\phi_i}([B_{(K_i)^\circ}]^\circ).  
\end{eqnarray*}  
If in addition all $\phi_i$ are increasing, then
\begin{eqnarray*}
\big[\OrliczAmix(\bK)\big]^{n} \leq  \big[\OrliczGmix(\bK)\big]^{n} \leq   \prod_{i=1}^n G^{orlicz}_{\phi_i}(B_{K_i})= \prod_{i=1}^n \O^{orlicz}_{\phi_i}(B_{K_i}).
\end{eqnarray*}   
 \et

In literature, the $i$-th mixed $p$-affine surface area and $i$-th mixed $L_p$ geominimal surface area are of interest, see details in e.g., \cite{Lut1987, WL1, WernerYe2010, DYzhuzhou2014, zhuzhou}. We will define and briefly discuss properties and inequalities for the Orlicz $i$-th mixed $L_{\phi}$ affine and geominimal surface areas. Let $K, L\in \cF_0^+$ and $Q_1, Q_2\in \cK_0$, we define $V_{\phi_1, \phi_2, i}(K, L; Q_1, Q_2)$ for $i\in \bbR$  by  
   \be \nonumber n V_{\phi_1, \phi_2, i}(K, L; Q_1, Q_2) \!=\!\! \int _{S^{n-1}}\!\! \bigg[\phi_1\!\!\left(\!\!\frac{h_{Q_1}(u)}{h_{K}(u)}\!\!\right)\! h_{K}(u)f_{K}(u)\bigg]^{\frac{n-i}{n}}\!\bigg[\phi_2\!\!\left(\!\frac{h_{Q_2}(u)}{h_{L}(u)}\!\!\right)\! h_{L}(u)f_{L}(u)\bigg]^{\frac{i}{n}} \!\! \,d\s(u).\ee For $L_1, L_2\in \cS_0$, we use the following formula   \be \nonumber n V_{\phi_1, \phi_2, i}(K, L; L_1^\circ, L_2^\circ) \!=\!\! \int _{S^{n-1}}\!\! \bigg[\phi_1\!\!\left(\!\!\frac{\r_{L_1}(u)^{-1}}{h_{K}(u)}\!\right)\! h_{K}(u)f_{K}(u)\bigg]^{\frac{n-i}{n}}\!\bigg[\phi_2\!\!\left(\!\frac{\r_{L_2}(u)^{-1}}{h_{L}(u)}\!\right)\! h_{L}(u)f_{L}(u)\bigg]^{\frac{i}{n}} \!\! \,d\s(u).\ee  H\"{o}lder's inequality (see
\cite{HLP}) implies 
\begin{eqnarray} \label{dual mixed volume-1}
[V_{\phi_1, \phi_2, i}(K, L; Q_1, Q_2)]^n &\leq&  [V_{\phi_1}(K, Q_1)]^{n-i}[V_{\phi_2}(L, Q_2)]^{i}, \ \  \ if
\ 0<i < n;\\
\nonumber [V_{\phi_1, \phi_2, i}(K, L; Q_1, Q_2)]^n &\geq&  [V_{\phi_1}(K, Q_1)]^{n-i}[V_{\phi_2}(L, Q_2)]^{i},\ \  \ if \ i <0 \ or \ i> n.\ \ \ \ \ 
\end{eqnarray} 
 
 We now propose our definition for the Orlicz $i$-th mixed $L_{\phi}$ affine surface area for $K, L\in \cF_0^+$.
\bd\label{mixed i th affine surface}  Let $K, L\in \cF_0^+$ and $i\in \bbR$. \\ 
(i)  For $\phi_1, \phi_2\in \Phi$, 
 \begin{eqnarray*} \O_{\phi_1, \phi_2,i}^{orlicz} (K, L) 
 &=&\inf_{ \{L_1, L_2\in \cS_0\} }
 \big\{nV_{\phi_1, \phi_2,i}\big(K, L; L_1^\circ, L_2^\circ\big): \ \ |L_1|=|L_2|=\o_n \big\}.
\end{eqnarray*}
(ii)   For $\phi_1, \phi_2\in \Psi$, 
 \begin{eqnarray*} \O_{\phi_1, \phi_2,i}^{orlicz} (K, L) 
 &=&\sup_{ \{L_1, L_2\in \cS_0\} }
 \big\{nV_{\phi_1, \phi_2,i}\big(K, L;  L_1^\circ,  L_2^\circ\big): \ \ |L_1|=|L_2|=\o_n \big\}.
\end{eqnarray*}
\ed    
Similarly,  the Orlicz $i$-th mixed $L_{\phi}$ geominimal surface area can be defined as follows.
\bd\label{mixed i th geominimal surface}  Let $K, L\in \cF_0^+$ and $i\in \bbR$.\\ 
(i)  For $\phi_1, \phi_2\in \Phi$, 
 \begin{eqnarray*} G_{\phi_1, \phi_2,i}^{orlicz} (K, L) 
=\inf_{ \{Q_1, Q_2\in \cK_0\} }
 \big\{nV_{\phi_1, \phi_2,i}\big(K, L; Q_1^\circ, Q_2^\circ\big): \ \ |Q_1|=|Q_2|=\o_n \big\}.
\end{eqnarray*}
(ii)   For $\phi_1, \phi_2\in \Psi$, 
 \begin{eqnarray*} G_{\phi_1, \phi_2,i}^{orlicz} (K, L) 
 =\sup_{ \{Q_1, Q_2\in \cK_0\} }
 \big\{nV_{\phi_1, \phi_2,i}\big(K, L;  Q_1^\circ,  Q_2^\circ\big): \ \ |Q_1|=|Q_2|=\o_n  \big\}.
\end{eqnarray*}
\ed    
 Clearly, $\O_{\phi_1, \phi_2,i}^{orlicz} (K, L)=\O_{\phi_1, \phi_2,n-i}^{orlicz} (L, K)$ and $G_{\phi_1, \phi_2,i}^{orlicz} (K, L)=G_{\phi_1, \phi_2,n-i}^{orlicz} (L, K)$ for all $i\in \bbR$ and all $K, L\in \cF_0^+$. Moreover,
\begin{eqnarray} \O_{\phi_1, \phi_2,0}^{orlicz} (K, L)=\O_{\phi_1}^{orlicz} (K)  \ \ \& \ \
 \O_{\phi_1, \phi_2, n}^{orlicz} (K, L)=\O_{\phi_2}^{orlicz} (L) ;\label{equation-3}\\ G_{\phi_1, \phi_2,0}^{orlicz} (K, L)=G_{\phi_1}^{orlicz} (K)  \ \ \& \ \
  G_{\phi_1, \phi_2, n}^{orlicz} (K, L)=G_{\phi_2}^{orlicz} (L) .\label{equation-4} \end{eqnarray}   
One can see that the Orlicz  $i$-th mixed $L_{\phi}$ affine and  geominimal surface areas  are all affine invariant. Moreover, for $K, L\in \cF_0^+$  and $i\in \bbR$ one has \begin{eqnarray} \O_{\phi_1, \phi_2,i}^{orlicz} (K, L)&\leq& G_{\phi_1, \phi_2, i}^{orlicz} (K, L), \ \ \phi_1, \phi_2\in \Phi;\label{compare:ith:----1}\\ \O_{\phi_1, \phi_2,i}^{orlicz} (K, L)&\geq& G_{\phi_1, \phi_2, i}^{orlicz} (K, L), \ \ \phi_1, \phi_2\in \Psi.\label{compare:ith:----2}  \end{eqnarray}   
 
\bt \label{cyclic-theorem} Let $K, L\in\cF_0^+$ and $i<j<k$. For $\phi_1, \phi_2\in\Psi$, one has 
 \begin{eqnarray*} \big[\O_{\phi_1, \phi_2,j}^{orlicz} (K, L)\big] ^{k-i}&\leq& \big[\O_{\phi_1, \phi_2,i}^{orlicz} (K, L)\big]^{k-j}\big[\O_{\phi_1, \phi_2,k}^{orlicz} (K, L)\big]^{j-i}; \\
  \big[G_{\phi_1, \phi_2,j}^{orlicz} (K, L)\big] ^{k-i}&\leq& \big[G_{\phi_1, \phi_2,i}^{orlicz} (K, L)\big]^{k-j}\big[G_{\phi_1, \phi_2,k}^{orlicz} (K, L)\big]^{j-i}. \end{eqnarray*}    \et

 \noindent \textbf{Proof.} We only prove the case for geominimal and omit the proof for affine case. Let $K, L\in \cF_0^+$ and $Q_1, Q_2\in\cK_0$. Let $i< j<k$ which implies
 $0<\frac{k-j}{k-i}<1$. Note also $k-i>0, k-j>0$ and $j-i>0$. H\"{o}lder's inequality  implies that,
 \begin{eqnarray}\nonumber
V_{\phi_1, \phi_2, j}(K, L; Q_1, Q_2)\leq [V_{\phi_1, \phi_2, i}(K, L; Q_1, Q_2)]^{\frac{k-j}{k-i}}[V_{\phi_1, \phi_2, k}(K, L; Q_1, Q_2)]^{\frac{j-i}{k-i}}.   
\end{eqnarray}
The desired result follows by taking the supremum over $Q_1, Q_2\in \cK_0$ with $|Q_1^\circ|=|Q_2^\circ|=\o_n$. 
 
\vskip 2mm \noindent{\bf Remark.}  Let  $\phi_1, \phi_2\in\Psi$. For $0< i< n$, let $(i, j, k) = (0, i, n)$ in Theorem \ref{cyclic-theorem},  by formulas (\ref{equation-3}) and (\ref{equation-4}), we have
 \begin{eqnarray*}
 \big[\O ^{orlicz}_{\phi_1, \phi_2, i}(K, L)\big]^{n}\!\leq\! \big[\O _{\phi_1}^{orlicz}(K)\big]^{n-i}\big[\O _{\phi_2}^{orlicz}(L)\big]^{i}  \ \&  \
\big[G ^{orlicz}_{\phi_1, \phi_2, i}(K, L)\big]^{n}\!\leq\! \big[G _{\phi_1}^{orlicz}(K)\big]^{n-i}\big[G _{\phi_2}^{orlicz}(L)\big]^{i}.
 \end{eqnarray*}
Equality always hold if  $i = 0$ or $i = n$. Similarly, for $i< 0$ or $i>n$, one has \begin{eqnarray}
\big[\O ^{orlicz}_{\phi_1, \phi_2, i}(K, L)\big]^{n}&\geq& \big[\O _{\phi_1}^{orlicz}(K)\big]^{n-i}\big[\O _{\phi_2}^{orlicz}(L)\big]^{i}, \nonumber \\ \big[G ^{orlicz}_{\phi_1, \phi_2, i}(K, L)\big]^{n}&\geq& \big[G _{\phi_1}^{orlicz}(K)\big]^{n-i}\big[G _{\phi_2}^{orlicz}(L)\big]^{i}.\label{Geominimal:ball:ith---1}
\end{eqnarray}
 
 \bt Let $K, L\in \cF_0^+$ be convex bodies with centroid (or Santal\'{o} point) at the origin.  \\ (i)  For $0\leq i \leq n$ and $\phi_1, \phi_2\in \Phi$,  one has  \begin{eqnarray*}\big[\O ^{orlicz}_{\phi_1, \phi_2, i}(K, L)\big]^{n}\leq \big[G ^{orlicz}_{\phi_1, \phi_2, i}(K, L)\big]^{n} \leq  \big[G_{\phi_1}^{orlicz}([B_{K^\circ}]^\circ)\big]^{n-i}\big[G_{\phi_2}^{orlicz}([B_{L^\circ}]^\circ)\big]^{i}.
   \end{eqnarray*}  
If in addition $\phi_1, \phi_2$ are both increasing, then \begin{eqnarray*}\big[\O ^{orlicz}_{\phi_1, \phi_2, i}(K, L)\big]^{n}\leq \big[G ^{orlicz}_{\phi_1, \phi_2, i}(K, L)\big]^{n} \leq  \big[G_{\phi_1}^{orlicz}(B_K)\big]^{n-i}\big[G_{\phi_2}^{orlicz}(B_L)\big]^{i}.
   \end{eqnarray*}  
 (ii)  Let $\E$ be an origin-symmetric ellipsoid. For $\phi_1, \phi_2\in \Psi$ and $i\leq 0$,   
\begin{eqnarray*} \big[\O_{\phi_1, \phi_2,i}^{orlicz} (K, \E)\big]^n &\geq& \big[G_{\phi_1, \phi_2, i}^{orlicz} (K, \E) \big]^n \geq   \big[G_{\phi_1}^{orlicz}(B_K) \big]^{n-i}\big[G _{\phi_2}^{orlicz}(\E)\big]^{i}. 
 \end{eqnarray*}  \et 
 \noindent  {\bf Proof.} (i). Let $\phi_1, \phi_2\in \Phi$  and $0\leq i\leq n$. Taking the infimum from both sides of inequality (\ref{dual mixed volume-1}) over $Q_1, Q_2\in\cK_0$ with $|Q_1^\circ|=|Q_2^\circ|=\o_n$ and by Definitions \ref{Orlicz geominimal surface} and  
\ref{mixed i th geominimal surface}, one has  \begin{eqnarray*}\big[G ^{orlicz}_{\phi_1, \phi_2, i}(K, L)\big]^{n}\leq \big[G _{\phi_1}^{orlicz}(K)\big]^{n-i}\big[G _{\phi_2}^{orlicz}(L)\big]^{i}.
 \end{eqnarray*}
 Combining with Theorem \ref{isoperimetric:geominimal} and inequality (\ref{compare:ith:----1}),  one gets, for $0\leq i \leq n$ and $\phi_1, \phi_2\in \Phi$,   \begin{eqnarray*}\big[\O ^{orlicz}_{\phi_1, \phi_2, i}(K, L)\big]^{n}  &\leq&  \big[G ^{orlicz}_{\phi_1, \phi_2, i}(K, L)\big]^{n}  \leq   \big[G _{\phi_1}^{orlicz}(K)\big]^{n-i}\big[G _{\phi_2}^{orlicz}(L)\big]^{i}\\ & \leq &  \big[G_{\phi_1}^{orlicz}([B_{K^\circ}]^\circ)\big]^{n-i}\big[G_{\phi_2}^{orlicz}([B_{L^\circ}]^\circ)\big]^{i}.
  \end{eqnarray*}  If in addition $\phi_1, \phi_2$ are both increasing, then \begin{eqnarray*}\big[\O ^{orlicz}_{\phi_1, \phi_2, i}(K, L)\big]^{n} &\leq&  \big[G ^{orlicz}_{\phi_1, \phi_2, i}(K, L)\big]^{n} \leq  \big[G _{\phi_1}^{orlicz}(K)\big]^{n-i}\big[G _{\phi_2}^{orlicz}(L)\big]^{i} \\ & \leq&   \big[G_{\phi_1}^{orlicz}(B_K)\big]^{n-i}\big[G_{\phi_2}^{orlicz}(B_L)\big]^{i}\!.
     \end{eqnarray*} 

\noindent (ii).    Let $i\leq 0$ and $\phi_1, \phi_2\in \Psi$. Note that both $\phi_1$ and $\phi_2$ are decreasing. 
 Inequalities (\ref{compare:ith:----2}) and (\ref{Geominimal:ball:ith---1}) together with Theorem \ref{isoperimetric:geominimal} imply 
 \begin{eqnarray*} \big[\O_{\phi_1, \phi_2,i}^{orlicz} (K, \E)\big]^n & \geq &  \big[G_{\phi_1, \phi_2, i}^{orlicz} (K, \E) \big]^n   \geq  \big[G _{\phi_1}^{orlicz}(K)\big]^{n-i}\big[G _{\phi_2}^{orlicz}(\E)\big]^{i} \\ &\geq&   \big[G_{\phi_1}^{orlicz}(B_K) \big]^{n-i}\big[G _{\phi_2}^{orlicz}(\E)\big]^{i}\!. 
    \end{eqnarray*}

\vskip 2mm \noindent {\bf Acknowledgments.} The research of DY is supported
 by a NSERC grant.

\vskip 2mm \noindent Deping Ye, \ \ \ {\small \tt deping.ye@mun.ca}\\
{\small \em Department of Mathematics and Statistics\\
   Memorial University of Newfoundland\\
   St. John's, Newfoundland, Canada A1C 5S7 }


\vskip 5mm \small

\begin{thebibliography}{M-PK.2}

\bibitem{Ale1937-1}
{A.D. Aleksandrov}, {\em On the theory of mixed volumes. i. Extension of certain concepts in
  the theory of convex bodies,} Mat. Sb. (N. S.) 2 (1937) 947-972. [Russian].

\bibitem{A1}
{   S. Alesker}, {\em Continuous rotation invariant valuations on
convex sets,} Ann. of  Math. 149 (1999) 977-1005.

\bibitem{A2}
{   S.\ Alesker}, {\em Description of translation invariant
valuations on convex sets with a solution of P. McMullen's
conjecture,} Geom. Funct. Anal. 11 (2001) 244-272.
 
\bibitem{Bl1}
{W. Blaschke}, {\em Vorlesungen {\"u}ber Differentialgeometrie II,
Affine Differentialgeometrie}, Springer Verlag, Berlin, 1923.

\bibitem{bor2012}
K.J. B\"{o}r\"{o}czky, E. Lutwak, D. Yang and G. Zhang, {\em The log-Brunn-Minkowski-inequality,} Adv. Math. 231 (2012) 1974-1997.

\bibitem{bor2013-1}
K.J. B\"{o}r\"{o}czky, {\em  Stronger versions of the Orlicz-Petty projection inequality,} J. Diff. Geom.  95 (2013)  215-247.


\bibitem{bor2013-2}
K.J. B\"{o}r\"{o}czky, E. Lutwak, D. Yang and G. Zhang, {\em The Logarithmic Minkowski Problem,} J. Amer. Math. Soc. 26 (2013)  831-852.


\bibitem{BM}
{J. Bourgain and V.D. Milman}, {\em New volume ratio
properties for convex symmetric bodies in $\bbR ^b$,} Invent.
Math. {88} (1987) 319-340.
 
 
 
 \bibitem{CG}
 {S. Campi and P. Gronchi,} {\em  The $L^p$ Busemann-Petty centroid
 inequality,} Adv. Math. 167 (2002) 128-141.
 
 \bibitem{Chen2011}
 {F. Chen, J. Zhou and C. Yang,} {\em On the reverse Orlicz
 Busemann-Petty centroid inequality,} Adv. Appl. Math. {47} (2011)
 820-828. 
 
  
 
\bibitem{Fenchel1938}
{W. Fenchel and B. Jessen}, {\em Mengenfunktionen und konvexe k\"{o}oper,} Danske Vid. Selskab. Mat.-fys. Medd. 16 (1938) 1-31.


\bibitem{Firey1962}
{W.J. Firey},  {\em $p$-means of convex bodies,}  Math. Scand. 10 (1962) 17-24.


\bibitem{Gardner2014}
 R.J. Gardner, D. Hug and W. Weil, {\em The Orlicz-Brunn-Minkowski theory: a general framework, additions, and inequalities,} J. Diff. Geom. in press. 
 
 
 \bibitem{Gr2}
 {  P.M. Gruber}, {\em Aspects of approximation of convex
 bodies,} Handbook of {C}onvex {G}eometry, vol. A, 321-345, North
 {H}olland, 1993.
 
 \bibitem{HaberlLYZ}
{C.\ Haberl, E.\ Lutwak, D.\ Yang and G.\ Zhang,} {\em The even
Orlicz Minkowski problem,} Adv.\ Math.\ {224} (2010) 2485-2510.


 \bibitem{HaberlFranz2009}
{C.\ Haberl and F.\ Schuster,} {\em General $L_p$ affine
isoperimetric inequalities,} J. Diff. Geom. {83} (2009) 1-26.

 

\bibitem{HLP}
{G.H. Hardy, J.E. Littlewood and G. P\'{o}lya}, {\em
Inequalities}, 2nd ed., Cambridge Univ. Press, 1952.



 \bibitem{Jenkinson2012}
{J. Jenkinson and E. Werner,} {\em Relative entropies for convex
bodies,}  to appear in Trans. Amer. Math. Soc.

\bibitem{GK2} 
{G. Kuperberg}, {\em From the Mahler conjecture to Gauss linking
integrals,} Geom. Funct. Anal. 18 (2008) 870-892.



 
 \bibitem{Li2011}
 {A. Li and G. Leng,} {\em A New Proof of the Orlicz Busemann-Petty
 Centroid Inequality,} Proc. Amer. Math. Soc. {139} (2011)
 1473-1481.
 
 
\bibitem{Ludwig2009}
{M. Ludwig,} {\em General affine surface areas}, Adv. Math. 224
(2010) 2346-2360.

\bibitem{LudR}
{M. Ludwig and M. Reitzner}, {\em A characterization of affine
surface area,} Adv. Math. 147 (1999) 138-172.

\bibitem{LR1}
{M. Ludwig and M. Reitzner,} {\em A classification of $SL(n)$
invariant valuations,} Ann. of  Math. 172 (2010) 1223-1271.



\bibitem {LSW}M. Ludwig, C. Sch\"{u}tt and E. Werner, \emph{Approximation of the Euclidean ball by polytopes}, Studia Math. 173 (2006) 1-18.

 
\bibitem{Lut1987}
{E. Lutwak}, {\em Mixed affine surface area,} J. Math. Anal. Appl.
125 (1987) 351-360.



\bibitem{Lut1993}
{E. Lutwak,} {\em The Brunn-Minkowski-Firey theory I: mixed volumes and the Minkowski problem,} J. Diff. Geom. 38 (1993) 131-150.


\bibitem{Lu1}{E. Lutwak}, {\em The Brunn-Minkowski-Firey theory. II. Affine and
geominimal surface areas,} Adv. Math. 118 (1996) 244-294.

\bibitem{LutwakZhang1997} 
E. Lutwak and G. Zhang, {\em Blaschke-Santal\'{o} inequalities}, J. Diff. Geom. 47 (1997) 1-16. 


  \bibitem{LYZ1}
 {E. Lutwak, D. Yang and G. Zhang,} {\em $L_p$ affine isoperimetric
 inequalities,} J. Diff. Geom.  56 (2000) 111-132.
 
 
 
 \bibitem{LYZ2010a}
 {E. Lutwak, D. Yang and G. Zhang,} {\em Orlicz projection bodies,}
 Adv. Math. {223} (2010) 220-242.
 
 \bibitem{LYZ2010b}
 {E. Lutwak, D. Yang and G. Zhang,} {\em Orlicz centroid bodies,}
 J. Diff. Geom. {84} (2010) 365-387.


\bibitem{MW2}
{M. Meyer and E. Werner}, {\em On the p-affine surface area,} Adv.
Math. 152 (2000) 288-313.

\bibitem{Nazarov2012}
 {F. Nazarov,}  {\em The H\"{o}rmander Proof of the Bourgain-Milman Theorem,}  
 Geom. Funct. Anal., Lecture Notes in Mathematics,  2050 (2012) 335-343.




\bibitem{Paouris2010}
{G. Paouris and E. Werner,} {\em Relative entropy of cone measures
and $L_p$ centroid bodies,} Proc. London Math. Soc. {104} (2012)
253-286.

\bibitem{Petty1974}
{C.M. Petty,} {\em Geominimal surface area,} Geom. Dedicata {3} (1974) 77-97.

\bibitem{Petty1985}
{C.M. Petty,} {\em Affine isoperimetric problems,} Annals of the New York Academy of Sciences, Volume 440, Discrete Geometry and Convexity, (1985) 113-127.

\bibitem{Sch} {R. Schneider,}  {\em Convex Bodies: The Brunn-Minkowski
theory,} Second edition, Cambridge Univ. Press, 2014.
 

\bibitem{Schneider2013}
 {R. Schneider,} {\em Affine surface area and convex bodies of elliptic type, } Periodica Math. Hungar.,
in press. 


 \bibitem{SW4}
{C. Sch{\"u}tt and E. Werner}, {\em Random polytopes of points
chosen from the boundary of a convex body,} in: GAFA Seminar
Notes, in: Lecture Notes in Math., vol. 1807, Springer-Verlag,
2002, pp. 241-422.


\bibitem{SW5}
{C. Sch{\"u}tt and E. Werner}, {\em Surface bodies and p-affine
surface area,} Adv. Math. 187 (2004) 98-145.

 
 \bibitem {WL1} W. Wang and G. Leng, \emph{$L_{p}$-mixed affine surface area}, J. Math. Anal. Appl. 335 (2007) 341-354.
 

\bibitem{Werner2012a}
{E. Werner,} {\em Renyi Divergence and $L_p$-affine surface area
for convex bodies,} Adv. Math. {230} (2012) 1040-1059.

\bibitem{Werner2012b}
{E. Werner,} {\em $f$-Divergence for convex bodies,} Proceedings of the ``Asymptotic Geometric Analysis" workshop, the Fields Institute, Toronto 2012.
 



\bibitem{WY2008}
{E. Werner and D. Ye}, {\em New $L_p$-affine isoperimetric
inequalities,} Adv. Math. {218} (2008) 762-780.

\bibitem{WernerYe2010}
{E. Werner and D. Ye}, {\em Inequalities for mixed $p$-affine
surface area,} Math. Ann. {347}   (2010) 703-737.


\bibitem {Y} D.\ Ye, \emph{Inequalities for general mixed affine surface areas}, J. London Math. Soc. 85 (2012) 101-120.

\bibitem{Ye2014} D. Ye, {\em $L_p$ Geominimal Surface Areas and their Inequalities,}  Int. Math. Res. Notes, in press. doi:10.1093/imrn/rnu009. 

\bibitem{Ye2013}
{D.\ Ye,} {\em On the monotone properties of general affine surface areas under the Steiner symmetrization,} Indiana Univ. Math. J., in press.  {\tt arXiv:1205.6145}


\bibitem{DYzhuzhou2014}
D. Ye, B. Zhu and J. Zhou, {\em On the mixed $L_p$ geominimal surface area for multiple convex bodies}, submitted. {\tt arXiv:1311.5180} 

\bibitem{Zhang2007}
{G. Zhang,} {\em New Affine Isoperimetric Inequalities,}
International Congress of Chinese Mathematicians (ICCM) {2}
(2007) 239-267.

 \bibitem{Zhou2011}
 {B. Zhu, N. Li and J. Zhou,} {\em Isoperimetric inequalities for $L_p$ Geominimal surface area,} Glasgow Math J. {53} (2011) 717-726.

\bibitem{zhuzhou} B. Zhu and J. Zhou, \emph{ $L_{p}$ mixed geominimal surface area}, submitted.

\bibitem{Zhu2012}
{G. Zhu}, {\em The Orlicz centroid inequality for star bodies,}
Adv. Appl. Math. {48} (2012) 432-445.


\end{thebibliography}
\end{document}